\documentclass[a4paper,12pt]{amsart}
\usepackage[a4paper,asymmetric]{geometry}
\usepackage[colorlinks]{hyperref}
\usepackage{mathrsfs}
\usepackage{graphicx}
\newtheorem{theorem}{Theorem}[section]
\newtheorem{lemma}[theorem]{Lemma}
\newtheorem{corollary}[theorem]{Corollary}
\newtheorem{proposition}[theorem]{Proposition}
\theoremstyle{remark}
\newtheorem{remark}[theorem]{Remark}
\theoremstyle{definition}
\newtheorem{definition}[theorem]{Definition}
\newtheorem{condition}[theorem]{Condition}
\numberwithin{equation}{section}

\DeclareMathOperator{\e}{e}
\newcommand{\spazio}{\mathcal{S}}
\newcommand{\map}{\mathcal{F}}
\newcommand{\Z}{\mathbf{Z}}
\newcommand{\N}{\mathbf{N}}
\newcommand{\R}{\mathbf{R}}

\newcommand{\scal}[1]{\langle #1\rangle}
\newcommand{\loc}{{\text{\Tiny loc}}}
\newcommand{\regular}{\mathscr{R}}
\newcommand{\singular}{\mathscr{S}}

\hypersetup{%
  pdftitle={Regularity and blow-up in a surface growth model},
  pdfsubject={The paper contains several regularity results and blow-up criterions
  for a surface growth model, which seems to have similar properties
  to the 3D Navier-Stokes, although it is a scalar equation.
  As a starting point we focus on energy methods and Lyapunov-functionals.},
  pdfauthor={D. Bloemker, M. Romito},
  pdfkeywords={surface growth, critical space, uniqueness, regularity, blow up, Leray estimates, Lyapunov function},
  pdfcreator={M. Romito}
}
\begin{document}
  \title[Surface growth blow up]{Regularity and blow-up in a surface growth model}
  \author[D. Bl\"omker]{Dirk Bl\"omker}
    \address{Institut f\"ur Mathematik\\ Universit\"at Augsburg\\ D-86135 Augsburg, Germany}
    \email{dirk.bloemker@math.uni-augsburg.de}
    \urladdr{http://www.math.uni-augsburg.de/ana/bloemker.html}
  \author[M. Romito]{Marco Romito}
    \address{Dipartimento di Matematica, Universit\`a di Firenze\\ Viale Morgagni 67/a\\ I-50134 Firenze, Italia}
    \email{romito@math.unifi.it}
    \urladdr{http://www.math.unifi.it/users/romito}
  \thanks{The second author gratefully acknowledges the support of \emph{Hausdorff Research Institute
    for Mathematics} (Bonn), through the \emph{Junior Trimester Program} on \emph{Computational Mathematics},
    and the hospitality of Augsburg Universit\"at.}
  \subjclass[2000]{35B33, 35B45, 35B65, 35K55, 35Qxx, 60H15}
  \keywords{surface growth, critical space, uniqueness, regularity, blow up, Leray estimates, Lyapunov function}
  \date{}
  \begin{abstract}
  The paper contains several regularity results and blow-up criterions
  for a surface growth model, which seems to have similar properties
  to the 3D Navier-Stokes, although it is a scalar equation.
  As a starting point we focus on energy methods and Lyapunov-functionals.
  \end{abstract}
  \maketitle
%%
%%%%%%%%%%%%%%%%%%%%%%%%%%%%%%%%%%%%%%%%%%%%%%%%%%%%%%%%%%%%%%%%%%%%%%%%%%%%%%%%%%%%%%%%%%%%%%%%%%%%%%%%%%%%%%%%%%%%

\tableofcontents

%%
%%
%%
%%%%%%%%%%%%%%%%%%%%%%%%%%%%%%%%%%%%%%%%%%%%%%%%%%%%%%%%%%%%%%%%%%%%%%%%%%%%%%%%%%%%%%%%%%%%%%%%%%%%%%%%%%%%%%%%%%%%
\section{Introduction}

Throughout this paper we consider a possible blow up for a model from surface
growth. Our main motivation is to carry over the program developed for 
3D-Navier stokes to this equation, in order to study the possible blow up
of solutions. This paper is the starting point focusing mainly on Hilbert space
theory.

Details on the model can be found in Raible et al.~\cite{Ra-Ma-Li-Mo-Ha-Sa:00},
\cite{Ra-Li-Ha:00} or Siegert \& Plischke \cite{Si-Pl:94:2nd}.
In its simplest version, it is given by
\begin{equation}\label{e:SG}
\partial_t h=-\partial_x^4 h -\partial_x^2 (\partial_x h)^2
\end{equation}
subject to periodic boundary conditions on $[0,L]$ and $\int_0^L h dx=0$.
Although the surface is not periodic, 
these boundary conditions together with the assumption of a moving frame
are the standard conditions in models of this type.
Sometimes the model has been considered also
on the whole real line without decay condition at infinity,
even though we do not examine this case here.

From a mathematical point of view Neumann
or Dirichlet boundary conditions are quite similar for the
problem studied here. The key point ensured by any of these
boundary conditions is that there is a suitable cancellation
in the non-linearity, namely
\begin{equation}\label{e:key}
 \int_0^L h\,(h_x^2)_{xx}\,dx =0\;,
\end{equation}
which is the main (and probably only) ingredient to 
derive useful a-priori estimates.

The main terms in the equation are the dominant linear operator, and the
quadratic non-linearity. Sometimes the equation is considered with 
a linear instability $-h_{xx}$,
which leads to the formation of hills, 
and the Kuramoto-Shivashinky-type  nonlinearity $(h_x)^2$
leading to a saturation in the coarsening of hills.
Both terms are  neglected here.
They are lower order terms not important for questions regarding regularity and blow up.
Moreover, the presence of these terms complicates calculations significantly
(cf. \cite{BlGuRa02}). 

Furthermore, the equation is usually perturbed 
by space-time white noise (see for instance~\cite{BlFlRo08}), which we also neglect here, although many results do hold for the stochastic PDE also.

For general surveys on surface growth processes and molecular beam
epitaxy see Barab\'asi \& Stanley \cite{BaSt95} or Halpin-Healy \& Zhang
\cite{HaZh95}.

%%%%%%%%%%%%%%%%%%%%%%%%%%%%%%%%%%%%%%%%%%%%%%%%%

\subsection{Existence of solutions}

%%%%%%%%%%%%%%%%%%%%%%%%%%%%%%%%%%%%%%%%%%%%%%%%%

There are two standard ways of treating the existence of solutions.
The first one relies on the spectral Galerkin method and shows energy 
type estimates for the approximation, which 
by some compactness arguments ensure the convergence of a subsequence.
See \cite{StWiP}, or for the stochastically 
perturbed equation \cite{BlGu02,BlHa04,BlFlRo08}.
In all cases initial conditions in $L^2$ ensure the existence, but not uniqueness, of
global solutions.

The second way uses fixed point arguments 
to show local uniqueness and regularity  using  the mild formulation.
See \cite{BlGu04}, which could not treat the optimal case.
In Section \ref{s:critical} we give a local existence,
which is optimal in the sense that initial conditions are 
in a critical space.
We also establish uniqueness among mild solutions and, less
trivially, among weak solutions. 
For these smooth local solutions we can easily
show energy estimates,
and discuss possible singularities and blow-up.

Standard arguments assure uniqueness of global solutions 
using a fixed point argument in $C^0([0,T],H^1)$ 
for sufficiently small regular data in $H^1$.
We can even go below that for uniqueness of  solutions
in $H^\alpha$ for any  $\alpha\geq\tfrac12$. This improves
results of~\cite{BlGu04}. But we are still not able to
prove uniqueness of global solutions without smallness
condition on the initial data. Nevertheless, we can give 
easily several conditions that imply uniqueness of global 
solutions. All of them assume regularity in critical spaces 
or more regularity (cf.\ Section \ref{sec:reg}).

In Section \ref{sec:blowup} we study possible 
singularities and blow up. Based on energy-type estimates,
we establish Leray-type estimates for lower bound on blow-up
in terms of $H^\alpha$-norms.
Moreover, we study an upper bound on the Hausdorff-dimension 
set of singularities in time, 
and show that a blow-up to $-\infty$ is more likely.
\begin{remark}
All results for regularity and Leray-type estimates 
 are based on energy estimates.
These are optimal in the sense that they hold also hold for complex
valued solutions. Furthermore, using the ideas of \cite{LiSin08},
\cite{LiSin08a}, one should be able to construct a complex valued
solution with strictly positive Fourier coefficients that actually
blows up in finite time. This is the subject of a work in progress.

This would show that results based on energy-estimates are useful
to describe a possible blow-up, but they alone will never be able
to rule it out.
\end{remark}
%%
%%%%%%%%%%%%%%%%%%%%%%%%%%%%
\subsection{Energy inequality}
\label{ss:ee}

We outline the standard idea for energy estimates,
which is to our knowledge the only useful idea for this equation.
If we formally multiply the equation by $h$ 
and integrate with respect to $x$, then we obtain using (\ref{e:key}),
\begin{equation}\label{e:energy} 
|h(t)|_{L^2}^2 + 2\int_0^t|\partial_x^2 h(s)|^2_{L^2}\,ds
\leq |h(0)|_{L^2}^2.
\end{equation}
Thus, using Poincare inequality,
$$
|h(t)|_{L^2} \le e^{-ct}  |h(0)|^2_{L^2}
\quad\text{and} \quad 
\int_0^\infty |h(t)|^2_{H^2}\,dt\leq  |h(0)|^2_{L^2}.
$$
As explained before this estimate is only valid for smooth local solutions,
or one could use spectral Galerkin approximation to verify it for global solutions.
Note that this regularity is lower than critical regularity.
It is enough for existence of solutions,
but not sufficient for uniqueness.
%%
%%
%%%%%%%%%%%%%%%%%%%%%%%%%%%%%%%%%%%%%%%%
\subsection{A Lyapunov-type functional}\label{sec:LF}

We can prove another a-priori estimate either for smooth local solutions or via 
spectral Galerkin approximations,
\begin{align*}
\tfrac1{\alpha^2}\partial_t \int_0^L \e^{\alpha h}\,dx
& = \int_0^L \e^{\alpha h}h_x h_{xxx}\,dx + 2\int_0^L \e^{\alpha h}h_x^2 h_{xx}\,dx\\
& = - \int_0^L \e^{\alpha h} h_{xx}^2 \,dx + (2 - \alpha) \int_0^L \e^{\alpha h}h_x^2  h_{xx}\,dx\\
& = - \int_0^L \e^{\alpha h} h_{xx}^2 \,dx - \tfrac13(2 - \alpha)\alpha  \int_0^L \e^{\alpha h}h_x^4\,dx.
\end{align*}
Thus, for $\alpha\in(0,2)$,
$$
\int_0^L \e^{\alpha h(t)}\,dx \leq  \int_0^L \e^{\alpha h(0)}\,dx \quad\text{for all }t>0
$$
and
$$
  \tfrac{(2 - \alpha)}{3}\alpha^3\int_0^\infty\int_0^L \e^{\alpha h}h_x^4\,dx\,dt
+ \alpha^2\int_0^\infty\int_0^L\e^{\alpha h} h_{xx}^2 \,dx\,dt
\leq \int_0^L \e^{\alpha h(0)}\,dx.
$$
With some more effort (cf. Stein-Winkler \cite{StWiP}), one knows that these terms are 
bounded independently  of $h(0)$ for large $t$.

The positive part $h^+=\max\{0,h\}$ now has much more regularity 
than the negative part $h^-=\max\{0,h\}$, so a possible blow up
seems to be more likely to $-\infty$ than to $+\infty$. 
We will illustrate this in Subsection \ref{ss:negblow}. But unfortunately,
this is still not sufficient regularity for uniqueness of solutions.

\begin{figure}
\includegraphics[scale=.5]{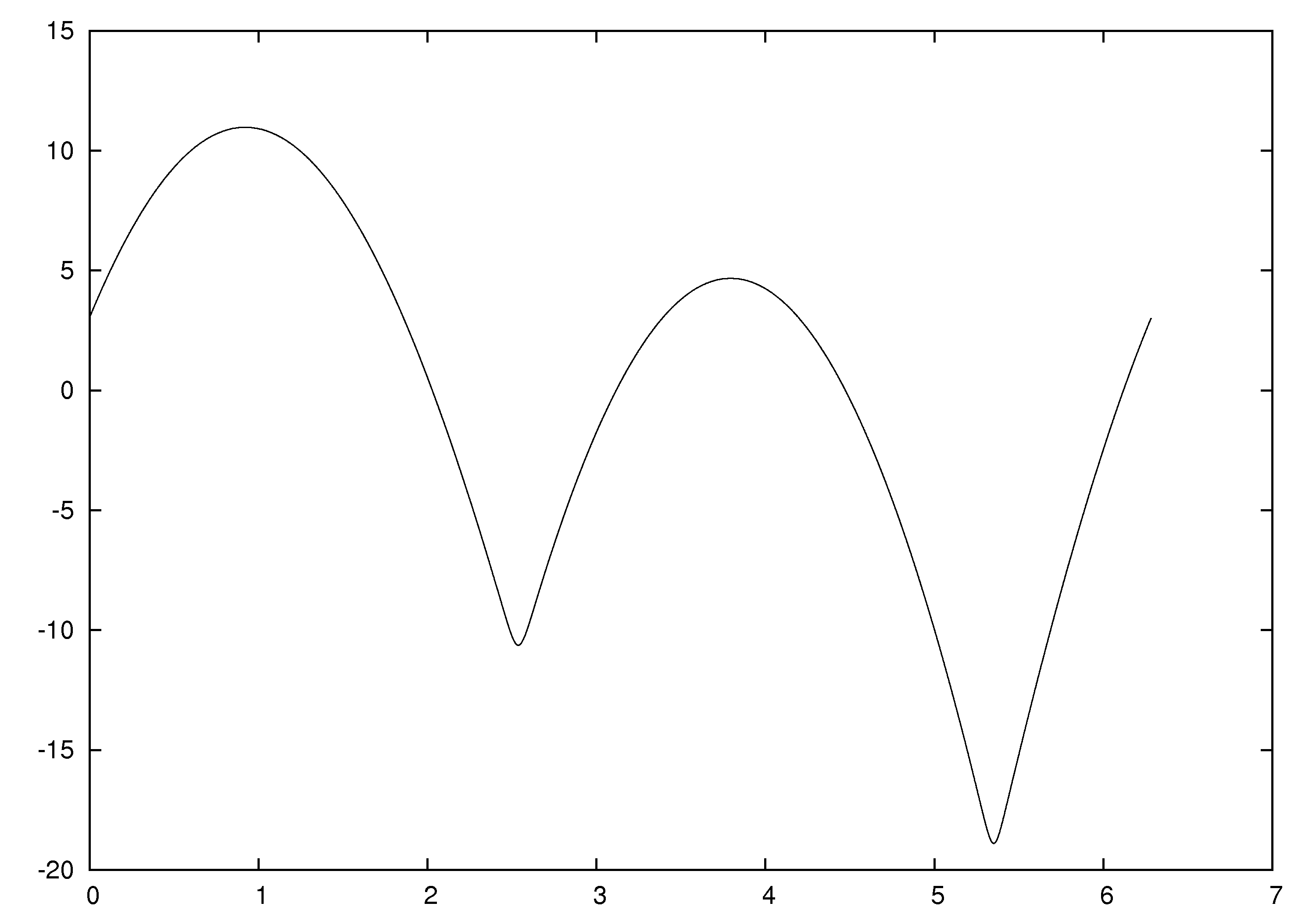}
\caption{A snapshot of a numerical solution 
to the surface growth equation 
with additional linear instability $-30\partial_x^2h$.
The hills look like parabola with 
sharp valleys in between.}
\end{figure}
%%
%%
%%
%%
%%%%%%%%%%%%%%%%%%%%%%%%%%%%%%%%%%%%%%%%%%%%%%%%%%%%%%%%%%%%%%%%%%%%%%%%%%%%%%%%%%%%%%%%%%%%%%%%%%%%%%%%%%%%%%%%%%%%
\section{Existence and uniqueness in a critical space}
\label{s:critical}

Prior to the details on some regularity criteria for equation \eqref{e:SG},
we introduce the \emph{scaling heuristic} which explains the formulae that
relate the different exponents in the results of the paper. An account on the
scaling heuristic for the Navier-Stokes equations can be found for
example in Cannone \cite{Can04}, such argument are on the ground
of the celebrated result on partial regularity for Navier-Stokes
of Caffarelli, Kohn \& Nirenberg \cite{CafKohNir82}. A recent paper
by Tao \cite{Tao08} discusses the scaling heuristic in the
framework of dispersive PDE.

The rationale behind the method is the following. First, notice
that the equations are invariant for the scaling transformation
\begin{equation}\label{e:scaling}
h(t,x) \longrightarrow h_\lambda(t,x)=h(\lambda^4 t, \lambda x).
\end{equation}
If $X$ is a functional space for $h$ (for example $L^\infty(0,T;L^2(0,L))$),
we can consider how the norm of $X$ scales with respect to the transformation
\eqref{e:scaling} above. Say the following relation holds,
$$
\|h_\lambda\|_X = \lambda^{-\alpha}\|h\|_X.
$$
We have the three cases
\begin{enumerate}
\item\emph{sub-critical} case for $\alpha<0$,
\item\emph{critical} case for $\alpha=0$,
\item\emph{super-critical} case for $\alpha>0$.
\end{enumerate}
The super-critical case corresponds to small-scales behaviour and is related
to low regularity, typically to topologies where possibly existence can be
proved, but no regularity or uniqueness. For example, one gets $\alpha=\frac12$
(hence, super-critical) for $X=L^\infty(0,\infty;L^2)$ or $X=L^2(0,T;\dot H^2)$,
which are the spaces where existence of global weak solutions can be proved.

The general scheme is the following. Consider spaces $X$ (depending on
the space variable) and $Y_T$ (depending on both variables, with $t$
up to $T>0$), then in order to have a regularity criterion based
on $Y_T$, the following statements must hold,
\begin{enumerate}
\item there is a unique local solution for every initial condition in $X$,
\item the unique local solution provided by (1) is regular,
\item the solution from (1) can be continued up to time $T$, as long
as its norm in $Y_T$ stays bounded.
\end{enumerate}

The above analysis has been extensively carried on by a large number
of authors for the three dimensional Navier-Stokes equations (see for
examples references in Cannone~\cite{Can04}). The first paper dealing
with such aims were Prodi~\cite{Pro59} and Serrin~\cite{Ser62}, see also
Beale, Kato \& Majda~\cite{BeaKatMaj84}.

%%
%%%%%%%%%%%%%%%%%%%%%%%%%%%%%%%
\subsubsection{Function spaces}

We shall mainly work in the hierarchy of Sobolev spaces of Hilbert type.
Since the equations are considered on $[0,L]$ with periodic boundary
conditions and zero space average, we shall use the following homogeneous
fractional Sobolev spaces. For $\alpha>0$,
$$
\dot H^{\alpha}
=\Bigl\{u\in L^2(0,L): u(\cdot+L)=u(\cdot),\quad u_0=0,\quad\sum_{k\not=0} k^{2\alpha}|u_k|^2<\infty\Bigr\},
$$
where $u_k$ is the $k^\text{th}$ Fourier coefficient, and $\dot H^{-\alpha}=(\dot H^\alpha)'$.
We shall consider the norm on $\dot H^\alpha$ defined by
\begin{equation}\label{e:normdotH}
|u|_{\alpha}^2 = \sum_{k\not=0} k^{2\alpha}|u_k|^2,
\end{equation}
which is equivalent to the norm of the Sobolev space $H^\alpha(0,L)$ on $\dot H^\alpha$.

We also use the space $L^p$ with norm $|\cdot|_{L^p}$ 
for the Lebesgue space of functions with integrable $p$-th power,
the space  $W^{k,p}$ with norm $|\cdot|_{W^{k,p}}$ for the 
Sobolev space, where the $k$-th derivative is in $L^p$, and the 
space $C^k$ of $k$-time continuously differentiable functions
with the supremum-norm.
%%
%%%%%%%%%%%%%%%%%%%%%%%%%%%%%%%%%%%%%%%%%%%%%%%%%%%%%%%%%%%%%%%%%%%%%%%%%%%%%%%%%%
\subsection{Existence and uniqueness in \texorpdfstring{$\dot H^{\frac12}$}{H1/2}}

This section is devoted to the proof of existence and uniqueness in the
critical space $\dot H^{\frac12}$, 
which improves significantly some results of Bl\"omker \& Gugg \cite{BlGu04}.
Here we shall follow the results of Fujita \& Kato~\cite{FujKat64}
on the Navier-Stokes equations with initial conditions in the
critical Sobolev Hilbert space. This is optimal in the sense
that local existence and uniqueness with lower regularity should imply uniqueness by rescaling.
\begin{definition}\label{d:Sa}
Given $T>0$, $\delta$ and $\alpha\in(0,\frac12)$, define the complete
metric space $\spazio_\alpha=\spazio_\alpha(T)$ as
$$
\spazio_\alpha(T)=\Big\{u\in C((0,T]; \dot H^{1+\alpha}):
\sup_{s\in(0,T]}\{s^{\frac{2\alpha+1}{8}}|u(s)|_{1+\alpha}\}<\infty\Big\},
$$
with norm
$$
\|u\|_{\alpha,T} =  \sup_{s\in(0,T]}\bigl\{s^{\frac{2\alpha+1}{8}}|u(s)|_{1+\alpha}\bigr\}
$$
and the $\delta$-ball 
$$
\spazio_\alpha^\delta(T) = \Big\{u\in\spazio_\alpha(T)\ :\  \|u\|_{\alpha,T}\leq \delta\Big\}.
$$
\end{definition}
Let us remark that for any $h\in \spazio_\alpha(T)$, $\widetilde\alpha\in(0,\alpha)$ and $\delta>0$
 we find  $\widetilde{T}\in(0,T)$ such that $h\in\spazio_{\widetilde\alpha}^\delta(\widetilde{T})$.
\begin{theorem}\label{t:criticalexuniq}
Given an arbitrary initial condition $h_0\in\dot H^{\frac12}$, there
exists a time $T_\bullet>0$, depending only on $h_0$, such
that there is a  solution $h\in C([0,T_\bullet);\dot H^{\frac12})$ to problem~\eqref{e:SG}.
Moreover,
\begin{enumerate}
\item $h\in C^\infty((0,T_\bullet)\times[0,L])$,
\item the solution satisfies the energy equality
      $$
      |h(t)|_{L^2}^2 + 2\int_0^t|h_{xx}|_{L^2}^2 = |h(0)|_{L^2}^2,
      $$
      for all $t<T_\bullet$,
\item there exists $a_\bullet>0$ such that $T_\bullet=+\infty$ if $|h_0|_{\frac12}\leq a_\bullet$.
\item Either the solution blows up in $\dot H^\beta$ for all $\beta>\frac12$ at $T=T_\bullet$  or $T_\bullet=\infty$.
\end{enumerate}
\end{theorem}
\begin{remark}\label{r:not_in_Hahalf}
If the maximal time $T_\bullet$ of a solution $h$ is finite, while we
know that $\|h(t)\|_\beta\to\infty$ as $t\uparrow T_\bullet$ for $\beta>\frac12$,
we cannot conclude that the same is true for $\|h(t)\|_{\frac12}$. Indeed,
$h$ can be discontinuous in the maximal time $T_\bullet$, so either
$\|h(t)\|_{\frac12}$ is unbounded, or is bounded and discontinuous in $T_\bullet$.

The reason behind this is that a solution in $\dot H^{\frac12}$ can be continued
as long as there is a control on the quantity $K_0$ of the type~\eqref{e:K0to0},
and this quantity is not uniformly convergent to $0$ in bounded subsets of $\dot H^{\frac12}$.
In different words, $K_0$ can be controlled as long as one can control the way
the mass of $h(0)$ is partitioned among Fourier modes.
\end{remark}
The proof of this theorem is developed in several steps, which we will prove in
the remainder of this section.

First, we prove existence and uniqueness (together with the global existence
statement). Then we prove an analogous result in $\dot H^\beta$, for all $\beta>\frac12$.
By a standard bootstrap technique, this implies the smoothness of solutions.

Let $A$ be the operator $\partial_x^4$ with domain $\dot H^4$. It is a standard
result that $A$ generates an analytic semigroup. Using
for example the Fourier series expansion,
it is easy to verify that
\begin{equation}\label{e:sgestimate}
|A^\gamma\e^{-tA}|_{\mathcal{L}(\dot{H}^\beta)}
\leq c_\gamma t^{-\gamma},
\end{equation}
for every $t>0$, where $\gamma\geq0$ and $\beta\in\mathbb{R}$. 
Moreover, it is easy to verify that the norm $|A^\frac{\beta}{4}\cdot|_{L^2} $,
which we will use several times in the paper, coincides with the standard
norm~\eqref{e:normdotH} on $\dot H^\beta$.

Proposition~\eqref{p:inequality} implies that for $\alpha\in(0,\frac12)$,
\begin{equation}\label{e:nonlin}
|A^{\frac18(4\alpha-5)}(h_x^2)_{xx}|_{L^2}
\leq c_\alpha |h|_{1+\alpha}^2
\end{equation}
(just apply the proposition with $\alpha=\beta$, $\gamma=\tfrac12-2\alpha$
and use the dual formulation of $L^2$ norm).

Consider now the right hand side of the mild formulation,
\begin{equation}\label{e:mild}
\map(h)(t)=\e^{-tA}h_0 + \int_0^t\e^{-(t-s)A}(h_x^2)_{xx}(s)\,ds,
\end{equation}
and define
\begin{gather*}
K_0(t)=\sup_{s\in(0,t]}\bigl(s^{\frac18(2\alpha+1)}|\e^{-sA}h_0|_{1+\alpha}\bigr),
	\quad\text{for }h_0\in \dot{H}^{\tfrac12}, \\
K(t,h)=\sup_{s\in(0,t]}\bigl(s^{\frac18(2\alpha+1)}|h(s)|_{1+\alpha}\bigr),
	\quad\text{for }h\in\spazio_\alpha(T),\ t\in[0,T].
\end{gather*}
Obviously, $K(t,h+k) \le K(t,h) + K(t,k) $
and
\begin{lemma}
For $h_0\in \dot{H}^{\tfrac12}$ we have
\begin{equation}\label{e:K0to0}
K_0(t)\to0\quad\text{as } t\to0.
\end{equation}
Furthermore, for each $\beta\in[\tfrac12,1+\alpha]$ there is a constant $c_\beta>0$ 
such that
\begin{equation}\label{e:K0bound}
K_0(t)\leq c_\beta t^{\frac18(2\beta-1)} |h_0|_\beta.
\end{equation}
\end{lemma}
\begin{proof}
By assumption $A^{\frac18}h_0\in L^2$, hence by Lemma~\ref{lem:SG} for $s\to 0$,
$$
s^{\frac18(2\alpha+1)}|\e^{-sA}h_0|_{1+\alpha}
=   |s^{\frac18(2\alpha+1)}A^{\frac18(1+2\alpha)}\e^{-sA} A^{\frac18}h_0|_{L^2} 
\to 0.
$$
For the second claim use (\ref{e:sgestimate}) to show
$$
K_0(t) 
=    \sup_{s\in(0,t]} s^{\frac18(2\alpha+1)}|A^{\frac14(1+\alpha-\beta)} e^{-sA}A^{\frac{\beta}{4}}h_0|_{L^2}
\leq c_\beta t^{\frac18(2\beta-1)}|h_0|_\beta.
$$
\end{proof}
Now we proceed to find a solution of $h=\map(h)$.
\begin{lemma}
There is a small constant $\delta>0$ depending on $\alpha$ such that for all
 $h_0\in \dot H^{\frac12}$ there exists a time $T$ sufficiently small,
such that the map $\map$ is a contraction on $\spazio_\alpha^\delta(T)$.
\end{lemma}
\begin{proof}
First we show that $\map$ maps $\spazio_\alpha^\delta$ into itself for $T$  and $\delta$ 
sufficiently small.
To be more precise, there is a number $c_\alpha>0$ such that for all $t\in[0,T]$
and all $h\in\spazio_\alpha^\delta$
\begin{equation}\label{e:claiminto}
K(t,\map(h))\leq K_0(t) + c_\alpha  K(t,h)^2 \leq K_0(T)+ c_\alpha \delta^2.
\end{equation}
Thus for $\delta\leq c_\alpha/2$ and $T$ sufficiently small $\map$ 
maps $\spazio_\alpha^\delta$ into itself.

In order to prove \eqref{e:claiminto} we consider
$$
|\map(h)(t)|_{1+\alpha}
\leq |\e^{-tA}h_0|_{1+\alpha} + \int_0^t|\e^{-(t-s)A}(h_x^2)_{xx}|_{1+\alpha}\,ds.
=I_0+I_1.
$$
For the first term,
$$
t^{\frac18(2\alpha+1)}I_0
=t^{\frac18(2\alpha+1)}|A^{\frac18(2\alpha+1)}\e^{-tA}h_0|_{\frac12}
\leq K_0(T)\to 0
$$
for $T\to0$.

For the second term we use \eqref{e:nonlin}, as well as \eqref{e:sgestimate}, to obtain
\begin{align*}
I_1
&=    \int_0^t|A^{\frac{5-4\alpha}8+\frac{1+\alpha}{4}}\e^{-(t-s)A}A^{\frac18(4\alpha-5)}(h_x^2)_{xx}|_{L^2}\,ds\\
&\leq c K(t, h)^2\int_0^t s^{-\frac14(2\alpha+1)}(t-s)^{-\frac18(7-2\alpha)}\,ds\\
&=C_\alpha t^{-\frac18(1+2\alpha)}K(t, h)^2,
\end{align*}
where $C_\alpha=cB(\frac14(3-2\alpha),\frac18(1+2\alpha))$ and
$B(x,y)=\int_0^1 t^{x-1}(1-t)^{y-1}\,dt$ is the \emph{Beta} function.

Now let us show that $\map$ is a contraction on $\spazio_\alpha$.
If  $h,k\in \spazio_\alpha$, then by following essentially the
above estimate of $I_1$, one can derive the following estimate
\begin{equation}\label{e:contraction}
K(T, \map(h)-\map(k))
\leq C_\alpha K(T, h-k)K(T, h+k) \leq 2\delta C_\alpha K(T, h-k)
\end{equation}
Thus $\map$ is a contraction, if $\delta\leq 1/(4C_\alpha)$.
\end{proof}

The following corollary is obvious, if we use \eqref{e:K0bound}
for $\beta>0$. The same conclusion cannot be drawn in the case
$\beta=\frac12$ (see Remark~\ref{r:not_in_Hahalf}).
\begin{corollary}
If $h_0\in \dot H^\beta$ for $\beta>\frac12$, then the time $T$
in the previous lemma depends only on a bound on $|h_0|_\beta$ and
not directly on $h_0$.
\end{corollary}
Thus, as long as a solution is bounded in any $\dot H^\beta$
with $\beta>\frac12$, the interval of existence can by extended
by a fixed length $T$, which depends only on the bounding constant.

The next lemma shows that the solution to the fixed point
$h=\map(h)$ in $\spazio_\alpha$ is continuous with values in $H^{1/2}$.
\begin{lemma}
If $h\in\spazio_\alpha(T)$, then $\map(h)\in C^0((0,T],\dot{H}^{1/2})$.
\end{lemma}
\begin{proof}
Obviously, it is enough to show that $\map(h)$ is continuous in $t=0$.
First, $e^{-tA}h_0\to 0$ in $\dot{H}^{\frac12}$ by continuity of the semigroup.
It remains to show that 
$$
\int_0^t e^{-(t-s)A}(h_x(s)^2)_{xx}\,ds \to 0
$$
in $\dot H^{\frac12}$ for $t\to0$. We know already by \eqref{e:nonlin} that 
$f(s) = s^{\frac14(2\alpha+1)}A^{\frac18(4\alpha-5)}(h_x^2)_{xx}$
is bounded in $L^2$ for $s\in(0,T]$ with $|f(s)|_{L^2}\leq cK(s,h)^2$. 
Thus from Lemma~\ref{lem:SG}),
$$
\int_0^ts^{-\frac14(2\alpha+1)} A^{-\frac18(4\alpha-5)+\frac18 } e^{-(t-s)A}f(s) \to 0
$$
in $L^2$, for $t\to0$.
\end{proof}
\begin{proposition}\label{p:critical}
Given $h_0\in\dot H^{\frac12}$ and $\alpha\in(0,\frac12)$, there exists
$T_0>0$ and $\delta_0$, depending only on $\alpha$ and $h_0$, such that
there is a unique solution in $\spazio_\alpha^{\delta_0}(T_0)$ to problem~\eqref{e:SG}
starting at $h_0$.

Moreover, the solution is in $C^0([0,T_0),\dot{H}^{1/2})$ and 
there exists $a_0>0$ small enough such that, if $|h_0|_{1/2}\leq a_0$,
then $T_0=\infty$.
\end{proposition}
\begin{proof}
Most of the proof is already done. We need to prove the last
statement of the proposition. By \eqref{e:sgestimate},
$K_0(t)\leq c_0 |h_0|_{\frac12}$, so that, if we choose
$a_0\leq (c_0c_\alpha)^{-1}$ (where $c_\alpha$ is the
constant in formula \eqref{e:claiminto}) and
$K = (2c_\alpha)^{-1}(1-\sqrt{1-c_0c_\alpha a_0})$,
by \eqref{e:claiminto} it follows that, for $K(t,h)\leq K$,
$$
K(t,\map(h))
\leq K_0(t) + c_\alpha K(t,h)^2
\leq c_0 a_0 + c_\alpha K^2
\leq K,
$$
independently of $t$. Hence, $T_0=\infty$.
\end{proof}
\begin{remark}[Criticality of $\spazio_\alpha(T)$]
Following the same notation used in Section \ref{s:critical},
we have that if $h\in\spazio_\alpha(T)$, then
$h_\lambda\in\spazio_\alpha(T_\lambda)$ and $K(T_\lambda,h_\lambda)$
scales as $\lambda^{\frac18(1-6\alpha)}K(T,h)$.
So, apparently, the $\|\cdot\|_{\alpha,T}$ does not obey
the scaling heuristic. On the other hand, this information
is of no use. Indeed, the scaling behaviour is hidden,
as it is shown by Lemma~\ref{l:addreg}, where the
boundedness in a space which is almost $\spazio_\alpha$
implies boundedness in the critical space
$L^q(0,T;\dot H^{1+\alpha})$, with $q=\frac8{1+2\alpha}$.
\end{remark}
Next, the case of more regular initial condition is considered.
The result is stated for integer exponents only, since for showing
regularity the present version  is sufficient (we already know
that solutions with initial value in $\dot H^{\frac12}$ are
continuous in $\dot H^1$). It is easy to adapt the proposition
to noninteger exponents, with some slight changes.
\begin{proposition}\label{p:critimorereg}
Let $n\in\N$, $n\geq1$. Given an arbitrary $h_0\in \dot H^n$,
there exist $T>0$ and a solution 
$h\in C([0,T);\dot H^n)\cap L^2_\text{loc}([0,T);\dot H^{n+2})$
to problem \eqref{e:SG}, with initial condition $h_0$.
\end{proposition}
\begin{proof}
We only prove the core \emph{a-priori} estimate for the
Theorem. Existence of a solution can be proven by means
of Proposition \ref{p:critical} or by an approximation
procedure (such as finite dimensional approximations).

Start by $n=1$,
$$
\frac{d}{dt}|h|_1^2
= 2\scal{h,\partial_t  h}_1
= - 2|h|_3^2 - 2\scal{h_{xx},(h_x^2)_{xx}}.
$$
By integration by parts and Sobolev,
interpolation and Young's inequalities, we get
$$
\begin{aligned}
2\scal{h_{xx},(h_x^2)_{xx}}
&=    - 2\scal{h_{xxx}, 2h_x h_{xx}}\\
	\text{\Tiny (by H\"older's inequality)}\quad
&\leq 2|h_{xxx}|_{L^2}|h_x|_{L^6}|h_{xx}|_{L^3}\\
	\text{\Tiny (by Sobolev embedding)}\quad
&\leq c|h|_3 |h|_{\frac43} |h|_{\frac{13}6}\\
	\text{\Tiny (by interpolation)}\quad
&\leq c|h|_3^{\frac74}|h|_1^{\frac54}\\
	\text{\Tiny (by Young's inequality)}\quad
&\leq |h|_3^2 + c|h|_1^{10}.
\end{aligned}
$$
In conclusion, if we denote by $\varphi(t)=|h(t)|_1^2+\int_0^t|h|_3^2$,
the above inequality reads
$$
\dot\varphi
=    \frac{d}{dt}|h|_1^2 + |h|_3^2
\leq c|h|_1^{10}
\leq c\varphi^5
$$
and by solving the differential inequality, we have a time $T$
such that $h$ is bounded in $C([0,T);\dot H^1)$ and in
$L^2_\text{loc}([0,T);\dot H^3)$.

The method is similar for $n\geq 2$. By computing the derivative
of $|h(t)|_n^2$, it turns out that it is necessary to estimate
the term originating from the nonlinear part. By integration
by parts and Leibnitz formula,
$$
\begin{aligned}
2\scal{D^{2n}h,(h_x^2)_{xx}}
&=    2\scal{D^{n+2}h,D^n(h_x^2)}\\
&=    2\sum_{k=0}^n\binom{n}{k}\scal{D^{n+2}h,(D^{k+1}h)(D^{n+1-k}h)}.
\end{aligned}
$$
By applying H\"older's inequality and Sobolev embedding, the above
sum can be estimated as above. All terms $|h|_a$ with $a\leq n$
can be controlled by $|h|_n$, while all terms with $a\in(n,n+2)$
can be controlled by $|h|_n$ and $|h|_{n+2}$ by interpolation. We
finally get the estimate
$$
\frac{d}{dt}|h|_n^2 + 2|h|_{n+2}^2
\leq |h|_{n+2}^2 + c_n |h|_n^{a_n},
$$
with suitable $c_n$ and $a_n$, depending only on $n$. By solving,
as above, the implied differential inequality, the solution $h$
turns out to be bounded in $C([0,T);\dot H^n)$ and in
$L^2_\text{loc}([0,T);\dot H^{n+2})$.
\end{proof}
Everything is now ready to carry on the proof of the main theorem
of this section.
\begin{proof}[Proof of Theorem \ref{t:criticalexuniq}]
The existence of solutions with initial condition in $\dot H^{1/2}$,
as well as the $T_\bullet=\infty$ statement, follow
from Proposition \ref{p:critical}.

The regularity statement (\textbf{1}) follows from Proposition
\ref{p:critimorereg}. Indeed, by Proposition \ref{p:critical},
a solution starting in $\dot H^{1/2}$ is continuous with
values in $\dot H^1$. By applying Proposition \ref{p:critimorereg}
on each $h(t)\in\dot H^1$, for $t\leq T_\bullet$, it follows that the solution
is $C((0,T_\bullet);\dot H^1)$ and $L^2_\text{loc}((0,T_\bullet);\dot H^3)$.
The last statement implies that $h(t)\in\dot H^3$, for almost every
$t\in(0,T_\bullet)$  and so Proposition \ref{p:critimorereg} can be
used with $n=3$, and so on. By iterating the procedure, it follows
that $h\in C((0,T_\bullet);\dot H^\beta)$ for all $\beta\geq1$.
Time regularity now follows from this space regularity and the
mild form~\eqref{e:mild}.

The energy equality in (\textbf{2}) is now easy using the space-time
regularity in $(0,T_\bullet)$ and the continuity at $t=0$ in the $L^2$ norm.
\end{proof}
%%
%%
%%
%%%%%%%%%%%%%%%%%%%%%%%%%%%%%%%%%%%%%%
\subsection{Uniqueness among weak solutions}

A weak solution to equation~\eqref{e:SG} is a function
$h\in L^\infty_\loc([0,\infty);L^2)\cap L^2_\loc([0,\infty);\dot H^2)$
which satisfies the equation in distributions. Existence of such
solutions for all initial data in $L^2$ has been established in
\cite{StWiP} (or \cite{BlGuRa02,BlFlRo08}). The following theorem shows that the solutions provided
by Theorem~\ref{t:criticalexuniq} are unique in the class of all weak
solutions $h$ that satisfy the energy inequality~\eqref{e:energy}.
\begin{theorem}\label{t:weakstrong}
Let $h_0\in\dot H^{\frac12}$ and let $h\in C([0,T_\bullet);\dot H^{\frac12})$ be the solution
to~\eqref{e:SG} provided by Theorem~\ref{t:criticalexuniq} and defined up to its maximal
time $T_\bullet$. Then every weak solution to~\eqref{e:SG} starting at $h(0)$ coincides with
$h$ on $[0,T_\bullet)$.
\end{theorem}
In order to prove the theorem, we shall proceed in several steps.
We will essentially prove that any solution in $\spazio_\alpha(T)$
with an additional integrability condition is unique in the class
of weak solutions (Proposition~\ref{p:weakstrong} below). Then we
prove that solutions in $\spazio_\alpha(T)$ satisfy the additional
condition (Lemma \ref{l:addreg} and \ref{l:last_step}). It is worth
remarking that the additional integrability condition~\eqref{e:ass_integ}
turns out to correspond to the critical space $L^{\frac{8}{1+2\alpha}}(H^{1+\alpha})$
(see Section~\ref{ss:criticality}).
\begin{proposition}\label{p:weakstrong}
Let $h\in \spazio_\alpha(T)$ be a solution to~\eqref{e:SG} and assume moreover that
\begin{equation}\label{e:ass_integ}
\int_0^T|h(t)|_{1+\alpha}^{\frac8{1+2\alpha}}\,dt
=\|h\|^{8/(2\alpha+1)}_{L^{8/(2\alpha+1)} ([0,T],H^{1+\alpha})} <\infty.
\end{equation}
Then $h$ is the unique weak solution starting at $h(0)$.
\end{proposition} 
\begin{proof}
Let $k$ be any weak solution starting at $h(0)$. Since $h\in C^\infty((0,T]\times[0,L])$
and $h$ is continuous in $\dot H^{\frac12}$, it follows that
$$
\scal{h(t),k(t)} + 2\int_0^t\scal{h_{xx},k_{xx}}\,ds
= -\int_0^t\int(h_{xx}k_x^2+k_{xx}h_x^2)\,ds
$$
which, together with the energy inequality for $k$ and the energy equality
(see Theorem~\ref{t:criticalexuniq}) for $h$ implies that the difference $w=h-k$
satisfies the following energy inequality,
$$
\begin{aligned}
|w(t)|_{L^2}^2 +2\int_0^t|w_{xx}|_{L^2}^2\,ds
&\leq 2\int_0^t\int (h_{xx}k_x^2+k_{xx}h_x^2)\,dx\,ds\\
&=    4\int_0^t\int k_xw_xw_{xx}\,dx\,dt
%&\leq |w_{xx}|_{L^2}^2 + c|w|_{L^2}^2|h|_{1+\alpha}^{\frac8{1+2\alpha}}.
\end{aligned}
$$
where we have used~\eqref{e:key} since
$$
h_{xx}k_x^2+k_{xx}h_x^2=2k_xw_xw_{xx}+w_{xx}w_x^2+h_{xx}h_x^2+k_{xx}k_x^2.
$$
The conclusion now follows from the assumption~\eqref{e:ass_integ} and
Gronwall's lemma, since
$$
4\int_0^t\int k_xw_xw_{xx}\,dx\,dt
\leq c|w|_2 |k|_{1+\alpha} |w|_{\frac{3-2\alpha}2}
\leq |w|_2^2 + c|k|_{1+\alpha}^{\frac{8}{1+2\alpha}}|w|_{L^2}^2,
$$
where we have used H\"older inequality (with exponents $2$, $\frac1{\alpha}$,
and $\frac{2}{1-2\alpha}$), the Sobolev embeddings $L^{\frac2{1-2\alpha}}\subset H^{\frac{3-2\alpha}2}$
and $L^{\frac1\alpha}\subset H^{1+\alpha}$, interpolation of $H^{\frac{3-2\alpha}2}$
between $L^2$ and $H^2$, and finally Young's inequality.
\end{proof}
Assumption~\eqref{e:ass_integ} cannot be obviously satisfied
by any arbitrary element of $\spazio_\alpha(T)$, hence we are led to
prove additional regularity for the solutions of~\eqref{e:mild}.
To this end, define for $T>0$ and $\alpha\in(0,\frac12)$,
$$
(\|u\|_{\alpha,T}^\star)^2
=\sum_{k\not=0} k^{2(1+\alpha)}\Big(\sup_{s\leq T} \{ s^{\frac18(1+2\alpha)}|u_k(s)|\}\Big)^2
$$
and
$$
\spazio_\alpha^\star(T) = \{u\in\spazio_\alpha(T): \|u\|_{\alpha,T}^\star<\infty\}.
$$
Assuming that $\spazio_\alpha^\star(T)\subset \spazio_\alpha(T)$ is not restrictive,
since it is easy to verify that $\|\cdot\|_{\alpha,T}\leq\|\cdot\|_{\alpha,T}^\star$.
\begin{lemma}\label{l:addreg}
If $h\in\spazio_\alpha^\star(T)$, then $\map(h)$ satisfies~\eqref{e:ass_integ} on $[0,T]$.
\end{lemma}
\begin{proof}
We write $\map(h)(t)=H_0(t)+H_1(t)$ where $H_0(t)=\e^{-tA}h(0)$
and $H_1$ contains the nonlinearity. Now,
$$
|H_0(t)|_{1+\alpha}^2 = \sum_{k\not=0} k^{2(1+\alpha)}\e^{-2ctk^4}|h_k(0)|^2,
$$
and so, if $\varphi\in L^q(0,T)$ with $p=\frac4{1+2\alpha}$ and $\frac1p+\frac1q=1$,
$$
\begin{aligned}
\int_0^T\varphi(t)|H_0(t)|_{1+\alpha}^2\,dt
&= \sum_{k\not=0}^\infty k^{2(1+\alpha)}|h_k(0)|^2\int_0^T\varphi(t)\e^{-2ctk^4}\,dt\\
&\leq \|\varphi\|_{L^q}\sum_{k\not=0}
 k^{2(1+\alpha)}|h_k(0)|^2\Bigl(\int_0^T\varphi(t)\e^{-2ctpk^4}\,dt\Bigr)^{\frac1p}\\
&\leq c_p\|\varphi\|_{L^q} |h(0)|_{\frac12}^2.
\end{aligned}
$$
By duality, the $L^{\frac8{1+2\alpha}}$ norm of $|H_0|_{1+\alpha}$ is finite. The second term is
more delicate, we shall proceed as in the proof of Proposition~\ref{p:inequality},
$$
\begin{aligned}
|H_1(t)|_{1+\alpha}^2
&= \sum_{k\not=0} k^{2(1+\alpha)}\Bigl(\int_0^t\e^{-c(t-s)k^4}[(h_x^2)_{xx}]_k\Bigr)^2\\
&= \sum_{k\not=0} k^{2(3+\alpha)}\Bigl(\sum_{l+m=k}|lm|\int_0^t\e^{-c(t-s)k^4}|h_l(s) h_m(s)|\,ds\Bigr)^2\\
&\leq \sum_{k\not=0} k^{2(3+\alpha)}\Bigl(\sum_{l+m=k}|lm|h_l^\star h_m^\star\Bigr)^2\Bigl(\int_0^t\e^{-c(t-s)k^4}s^{-\frac{1+2\alpha}4}\,ds\Bigr)^2,
\end{aligned}
$$
where $h_k^\star=\sup_{s\leq T}s^{\frac{1+2\alpha}8}|h_k(s)|$. Hence, for every
$\varphi\in L^q(0,T)$,
\begin{multline*}
\int_0^T\varphi(t)|H_1(t)|_{1+\alpha}\,dt\leq\\
\leq \sum_{k\not=0} k^{2(3+\alpha)}\Bigl(\sum_{l+m=k}|lm|h_l^\star h_m^\star\Bigr)^2
          \int_0^T\varphi(t)\Bigl(\int_0^t\e^{-c(t-s)k^4}s^{-\frac{1+2\alpha}4}\Bigr)^2\,dt.
\end{multline*}
If we prove that
\begin{equation}\label{e:elem_claim}
\int_0^T\varphi(t)\Bigl(\int_0^t\e^{-c(t-s)k^4}s^{-\frac{1+2\alpha}4}\,ds\Bigr)^2\,dt
\leq c\|\varphi\|_{L^q}k^{2\alpha-7},
\end{equation}
then we can proceed as in the proof of Proposition~\ref{p:inequality} (where the $h_k^\star$
replace the Fourier components and $\gamma=\frac12-2\alpha$) to obtain that
$$
\int_0^T\varphi(t)|H_1(t)|_{1+\alpha}\,dt
\leq c\|\varphi\|_{L^q}(\|h\|_{\alpha,T}^\star)^2,
$$
and, again by duality, boundedness of $\map(h)$.

So, everything boils down to proving~\eqref{e:elem_claim}. Using H\"older
inequality and (twice) a change of variables,
$$
\begin{aligned}
\lefteqn{\int_0^T\varphi(t)\Bigl(\int_0^t\e^{-c(t-s)k^4}s^{-\frac{1+2\alpha}4}\Bigr)^2\,dt\leq}\qquad\\
&\leq \|\varphi\|_{L^q}\Bigl[\int_0^T\varphi(t)\Bigl(\int_0^t\e^{-c(t-s)k^4}s^{-\frac{1+2\alpha}4}\,ds\Bigr)^{2p}\,dt\Bigr]^{\frac1p}\\
&\leq \|\varphi\|_{L^q} k^{2\alpha-7}\Bigl[\int_0^\infty\Bigl(\int_0^t\e^{-(t-s)}s^{-\frac{1+2\alpha}4}\,ds\Bigr)^{2p}\,dt\Bigr]^{\frac1p},
\end{aligned}
$$
and it is elementary to verify that the integral on the right-hand side is convergent.
Indeed,
$$
\int_0^{\frac{t}2}\e^{-(t-s)}s^{-\frac{1+2\alpha}4}\,ds
\leq c t^{\frac{3-2\alpha}4}\e^{-\frac{t}2},
$$
which is in $L^p(0,\infty)$, as well as
$$
\int_{\frac{t}2}^t\e^{-(t-s)}s^{-\frac{1+2\alpha}4}\,ds
\leq c t^{-\frac{2\alpha+1}4}(1-\e^{-\frac{t}2}),
$$
since $t^{-2p\frac{2\alpha+1}4}=t^{-2}$.
\end{proof}
The final step is to prove that solutions exist in the smaller space
$\spazio_\alpha^\star$. This is then the unique weak solution and the
solution given by Theorem \ref{t:criticalexuniq}.
\begin{lemma}\label{l:last_step}
Let $h_0\in\dot H^{\frac12}$ and $\alpha\in(0,\frac12)$. Then there
is $T_\star>0$ such that there exists a solution $h$ in $\spazio_\alpha^\star(T_\star)$.
\end{lemma}
\begin{proof}
The proof is essentially a fixed point argument, as in Proposition~\ref{p:critical}.
So, it is sufficient to show the following facts:
\begin{enumerate}
\item $\|H_0\|_{\alpha,T}^\star\leq\|h(0)\|_{\frac12}$,
\item $\|H_0\|_{\alpha,T}^\star\longrightarrow0$ as $T\to0$,
\item there is $c>0$ (independent of $T$) such that for all $h\in\spazio_\alpha^\star(T)$,
      $\|\map(h)\|_{\alpha,T}^\star\leq \|H_0\|_{\alpha,T}^\star + c(\|h\|_{\alpha,T}^\star)^2$,
\item there is $c>0$ (independent of $T$) such that $\|\map(g)-\map(h)\|_{\alpha,T}^\star\leq c\|g-h\|_{\alpha,T}^\star \|g+h\|_{\alpha,T}^\star$
      for all $g$, $h\in\spazio_\alpha^\star(T)$,
\end{enumerate}
where $H_0(t)=\e^{-tA}h(0)$ and $H_1=\map(h)-H_0$. Notice that
$$
\sup_{s\leq T}s^{\frac{1+2\alpha}8}|[H_0(t)]_k|
= |h_k(0)|\sup_{s\leq T}s^{\frac{1+2\alpha}8}\e^{-csk^4}
\leq c k^{-\frac{1+2\alpha}2}|h_k(0)|
$$
and so
$$
(\|H_0\|_{\alpha,T}^\star)^2
\leq \sum_{k\not=0} k^{2(1+\alpha)}|h_k(0)|^2 ck^{-(1+2\alpha)}
\leq c\|h(0)\|_{\frac12}^2.
$$
In order to prove the second property, we have to refine the previous
computation. Fix $\varepsilon>0$ such that $\varepsilon\leq c_\alpha$
(where $c_\alpha^4$ is the point where the function $s^{\frac{1+2\alpha}8}\e^{-s}$
attains its maximum), then 
$$
\begin{aligned}
(\|H_0\|_{\alpha,T}^\star)^2
&=   \Bigl(\sum_{|k|\leq\varepsilon T^{-\frac14}} + \sum_{|k|>\varepsilon T^{-\frac14}}\Bigr)
        k^{2(1+\alpha)}|h_k(0)|^2(\sup_{s\leq T} s^{\frac{1+2\alpha}4}\e^{-2csk^4})\\
&\leq  \sum_{|k|\leq\varepsilon T^{-\frac14}}(k^4T)^{\frac{1+2\alpha}4} |k||h_k(0)|^2
    + \sum_{|k|>\varepsilon T^{-\frac14}} k|h_k(0)|^2\\
&\leq \varepsilon^{1+2\alpha}\|h(0)\|_{\frac12}^2 + c\sum_{|k|>\varepsilon T^{-\frac14}} |k||h_k(0)|^2.
\end{aligned}
$$
Now, $\limsup_{T\to0}\|H_0\|_{\alpha,T}^\star\leq\varepsilon^{1+2\alpha}\|h(0)\|_{\frac12}^2$
and, as $\varepsilon\downarrow0$, the conclusion follows.

In order to prove the last fact, we follow the proof
of Lemma~\ref{l:addreg},
$$
|(H_1)_k(t)|
\leq k^2\sum_{l+m=k}|l m| h_l^\star h_m^\star\Bigl(\int_0^t\e^{-c(t-s)k^4}s^{-\frac{1+2\alpha}{4}}\Bigr)
$$
and so
$$
\sup_{t\leq T}t^{\frac{1+2\alpha}8}|(H_1)_k(t)|
\leq k^2\sum_{l+m=k}|l m| h_l^\star h_m^\star\Bigl(\sup_{t\leq T}t^{\frac{1+2\alpha}8}\int_0^t\e^{-c(t-s)k^4}s^{-\frac{1+2\alpha}{4}}\Bigr).
$$
Assume that the term in round brackets in the above formula is bounded by $c_\alpha k^{-\frac{7-2\alpha}2}$
(we shall prove this later), then, as in the proof of Proposition~\ref{p:inequality},
$$
(\|H_1\|_{\alpha,T}^\star)^2
\leq c_\alpha\sum_{k\not=0} k^{6+2\alpha}\Bigl(k^{-\frac{7-2\alpha}2}\sum_{l+m=k}|l m| h_l^\star h_m^\star\Bigr)^2
\leq c(\|h\|_{\alpha,T}^\star)^4.
$$
As it regards the rounded brackets term, we use the inequality $\e^{-c(t-s)k^4}\leq c_\alpha [k^4(t-s)]^{-\frac{7-2\alpha}{8}}$
to get
$$
\Bigl(\sup_{t\leq T}t^{\frac{1+2\alpha}8}\int_0^t\e^{-c(t-s)k^4}s^{-\frac{1+2\alpha}{4}}\,ds \Bigr)
\leq c_\alpha B(\tfrac{3-2\alpha}4,\tfrac{1+2\alpha}8)k^{-\frac{7-2\alpha}2},
$$
and $B$ is the \emph{Beta} function.

The proof of the last fact is similar. Indeed, if $g$, $h\in\spazio_\alpha^\star(T)$, then
$$
\begin{aligned}
|[\map(g)(t)-\map(h)(t)]_k|
&\leq k^2\int_0^t\e^{-c(t-s)k^4}|[(g-h)_x(g+h)_x]_k|\,ds\\
&\leq k^2\sum_{l+m=k}|l m|\int_0^t\e^{-c(t-s)k^4}|(g_l-h_l)(g_m+h_m)|\,ds
\end{aligned}
$$
and so, by proceeding as above, the last fact follows.
\end{proof}
\begin{proof}[Proof of Theorem~\ref{t:weakstrong}]
Given $h(0)\in\dot H^{\frac12}$, let $h\in C([0,T_\bullet);\dot H^{\frac12})$
be the solution provided by Theorem~\ref{t:criticalexuniq} and fix $T<T_\bullet$.
By Lemma~\ref{l:last_step} we know that $h\in\spazio_\alpha^\star(T_\star)$,
so Lemma~\ref{l:addreg} implies that $h$ satisfies the integrability
condition~\eqref{e:ass_integ} on $[0,T_\star]$. By property (\textbf{1})
of Theorem~\ref{t:criticalexuniq}, $h$ satisfies trivially~\eqref{e:ass_integ}
on $[T_\star,T]$. So Proposition~\ref{p:weakstrong} applies and the
conclusion follows.
\end{proof}
%%
%%
%%
%%%%%%%%%%%%%%%%%%%%%%%%%%%%%%%%%%%%%%%%%%%%%%%%%%%%%%%%%%%%%%%%%%%%%%%%%%%%%%%%%%%%%%%%%%%%%%%%%%%%%%%%%%%%%%%%%%%%
\section{Regularity}
\label{sec:reg}
%%
%%
%%
%%%%%%%%%%%%%%%%%%%%%%%%%%%%%%%%%%%%%%%%%%%%%%%%%%%%%%%%%%
\subsection{Criticality}\label{ss:criticality}

In this section, we carry out the program described in the beginning of the previous section. 
We will find spaces $Y_T$ such that boundedness in these spaces imply 
uniqueness for solutions starting in $H^{1/2}$.

Let us first discuss regularity criteria in Lebesgue spaces.
Set $T_\lambda=\lambda^{-4}T$ and $L_\lambda=\lambda^{-1}L$ and consider the space
$X(\lambda)=L^q(0,T_\lambda;L^p(0,L_\lambda))$, for some values of $p$ and $q$.
Under the scaling \eqref{e:scaling} we have that
$$
\|h_\lambda\|_{X(\lambda)} = \lambda^{-\frac{4}{q}-\frac{1}{p}}\|h\|_{X(1)}.
$$
so that the space $L^\infty((0,T)\times(0,L))$ turns out to be the only critical
space in this class. All other Lebesgue spaces are super-critical.

The conjecture now is that solutions in $L^\infty(0,T;L^\infty(0,L))$ 
or $C((0,T)\times(0,L))$ are unique and regular. 
We believe that with similar methods, as in the existence for 
initial conditions in $H^{1/2}$, one should be able to 
prove existence of unique local solutions. 
But this  is much more involved.

In order to consider Sobolev spaces,  we set
$X(\lambda)=L^q(0,T_\lambda;\dot W^{k,p}(0,L_\lambda))$ and
$$
\|h_\lambda\|_{X(\lambda)} = \lambda^{k-\frac{4}{q}-\frac{1}{p}}\|h\|_{X(1)}.
$$
(this is easy for integer $k$ and tricky for non-integer values, but
it can be done). Hence, the space is critical for
$$
\frac{4}{q} + \frac{1}{p} = k.
$$
In the following subsection, we will give 
the corresponding criteria for $p=2$, $k$ arbitrary
and $p=4$, $k=1$. 
The extension to $k=1$ and $p$ arbitrary is straightforward
and not presented here.

Let us finally remark, that in the following,
we also give regularity criteria for $L^4(0,T,C^1(0,L))$,
which is also a critical space. 
%
%%
%%%%%%%%%%%%%%%%%%%%%%%%%%%%%%%%%%%%%%%%%%%%%%%%%%%%%%%%%%
\subsection{Regularity Criteria}

In principle the following {\em Meta-theorem} should hold:
If a solution is bounded in a critical space, 
then it is unique, and does not have a blow up.
This means that the unique local solution exists as long
as at least one (hence all, as the solution is then proved to be regular)
of the critical norms is finite over the time horizon.

For simplicity, in the rest of the section we focus only
on some examples and we consider solutions with 
sufficiently smooth initial condition, in order to have energy type 
estimates for the $H^1$-norm without any trouble at $t=0$.

We just remark that energy estimates in any 
other $H^s$-space with $s>\tfrac12$ yield exactly the same result.

\begin{theorem}\label{t:regul}
Let $h_0\in\dot H^1$, let $h=h(\cdot,h_0)$ be the unique local
solution started at $h_0$ and let $\tau(h_0)$ be the maximal
time of $h$. Then $h$ is $C^\infty$ in
space and time on $(0,\tau(h_0))$ and for every $\alpha\in(\frac12,\frac92)$,
$$
\int_0^{\tau(h_0)}\|h(s)\|_{H^\alpha}^{\frac{8}{2\alpha-1}}\,ds
=\infty.
$$
Moreover, 
$$
\int_0^{\tau(h_0)}\|h(s)\|_{W^{1,4}}^{16/3}\,ds
=\infty
\quad\text{and}\quad
\int_0^{\tau(h_0)}\|h(s)\|_{C^1}^4\,ds
=\infty\;.
$$
\end{theorem}
\begin{proof}
We already know by Theorem \ref{t:criticalexuniq}
that there is a unique local solution in $C((0,\tau); \dot H^1)$
for initial conditions in $\dot H^1$, which is actually smooth. Indeed  
$h\in C^\infty((0,\tau)\times(0,L))$. Furthermore, the $H^1$-norm blows up at $t\to \tau$.

Now fix $\alpha\in(\frac12,\frac92)$, then by integration by parts and the Sobolev embedding
$H^\frac16 \subset L^3$,
\begin{align} \label{e:crit_H}
\frac{d}{dt}|h|_{1}^2 + 2|h|_{3}^2
& = -2\langle h, (h_x^2)_{xx}\rangle_{1}
  = -2\langle h_x, (h_x^2)_{xxx}\rangle_{L^2}\nonumber\\
& = -4\int_0^L h_x h_{xx} h_{xxx}\,dx
  = 2\int_0^L h_{xx}^3\nonumber\\
&\leq c|h|_{{\frac{13}6}}^3.
\end{align}
By interpolation, it is easy to see that
$$
|h|_{\frac{13}6}^3\leq |h|_{1}^{\frac{2\alpha-1}4}|h|_{\alpha}|h|_{3}^{\frac{9-2\alpha}4},
$$
and so using Young's inequality,
$$
\frac{d}{dt}|h|_{1}^2 + 2|h|_{3}^2
\leq |h|_{3}^2 + c|h|_{\alpha}^{\frac8{2\alpha-1}}|h|_{1}^2.
$$
Finally, by Gronwall's lemma, the proof of the first statement is complete.

Let us turn again to (\ref{e:crit_H}). 
Using Sobolev embedding $H^\frac14\subset L^4$ yields
\begin{equation}\label{e:crit_W}
 \frac{d}{dt}|h|_{1}^2 + 2|h|_{3}^2 \le C|h|_{W^{1,4}}|h|_{2+\frac14}|h|_3\;.
\end{equation}
Again by interpolation and Young inequality
$$
\frac{d}{dt}|h|_{1}^2 + 2|h|_{3}^2
\leq |h|_{3}^2 + C|h|_{W^{1,4}}^{16/3}|h|_{1}^2\;,
$$
which yields the result using Gronwall.

The last claim follows similarly, using
\begin{equation}\label{e:crit_C}
 \frac{d}{dt}|h|_{1}^2 + 2|h|_{3}^2 \le C|h|_{C^1}|h|_2|h|_3\;.
\end{equation}
\end{proof}
%
% 
%%
%%%%%%%%%%%%%%%%%%%%%%%%%%%%%%%%%%%%%%%%%%%%%%%%%%%%
\subsection{\texorpdfstring{$H^3$}{H3}-regularity}

In this section we show $L^p(0,T, H^3)$ 
for some small $p$ which is possibly  less than $1$.
We gain spatial regularity by paying time regularity.
The main result is:
\begin{theorem}
If for a solution $h\in L^r(0,T,H^1)$ for some $r\in(0,10)$, then  $h\in L^{r/5}(0,T,H^3)$.
Moreover,
$$
 \int_0^T |h_{xxx}|^{r/5} \,dt 
\le C \Big(\int_0^T |h_x|^r \, dt \Big)^{(10-r)/5}\;.
$$
\end{theorem}
\begin{remark}
It is easy to check that the space $L^{r/5}(0,T,H^3)$ is critical
if and only if $L^r(0,T,H^1)$ is critical. Thus this result respects
the criticality heuristic.
\end{remark}
\begin{remark}
If $h\in L^\infty(0,T,H^{1/2})$ (critical)
then by interpolation of  $H^{1/2}$ and  $H^2$ we obtain  from energy estimates $h\in L^6(0,T,H^1)$,
 and now  $h\in L^{6/5}(0,T,H^3)$.
Then by interpolation of  $H^{1/2}$ and $H^3$  we recover $h\in L^2(0,T,H^2)$.
Thus this regularity result gives no improvement of the regularity given by the 
energy estimate in Section \ref{ss:ee}. It respects the level 
of criticality of the spaces.
\end{remark}
\begin{proof}
For some $p>0$ where $|\cdot|$ denotes the norm in $L^2$
$$
\partial_t \frac{|h_x|^2}{1+|h_x|^p}
 =   2 \frac{\langle h_{xt},h_x\rangle}{1+|h_x|^p} + p \frac{|h_x|^2 \langle h_{xt},h_x \rangle}{(1+|h_x|^p)^2}\\
 =:   \varphi_p(|h_x|)  \langle h_{xt},h_x\rangle 
$$
where $ \varphi_p(z)=(2+(2+p)z^2)/(1+z^p)^2$. 
Thus using the PDE and integration by parts
$$
\partial_t \frac{|h_x|^2}{1+|h_x|^p} 
=  \varphi_p(|h_x|) (\tfrac12 \int h_{xx}^3\,dx   - |h_{xxx}|^2)
$$
Using the embedding of $H^{1/6}$ into $L^3$,  interpolation, and Young  yields 
$$
| \int h_{xx}^3\,dx| 
\le C |h|_{{13/6}}^3 
\le C |h|_{1}^{5/4} |h|_{3}^{7/4}
 \le  C |h|_{1}^{10} +  |h|_{3}^2 
$$
Combining both results yields 
$$
\partial_t \frac{|h_x|^2}{1+|h_x|^p} = -\tfrac12\varphi_p(|h_x|) |h_{xxx}|^2 + C \varphi_p(|h_x|) |h_x|^{10}
$$
Now, $\varphi_p(|h_x|) |h_x|^{10}\in L^1(0,T)$ if $12-r-2p\le0$ and thus $p\ge 6-\tfrac{r}2$.
We derive 
$$ \int_0^T \tfrac12\varphi_p(|h_x|) |h_{xxx}|^2 \,dt <\infty \quad\text{for}\quad p= 6-\tfrac{r}2>1.
$$
Using H\"older inequality for some $\alpha\in(0,2)$ yields
\begin{eqnarray*}
  \int_0^T |h_{xxx}|^\alpha \,dt 
&\le&\Big(\int_0^T |h_{xxx}|^2\varphi_p(|h_x|)  \,dt \Big)^{\alpha/2}
 \Big(\int_0^T\varphi_p(|h_x|)^{-\alpha/(2-\alpha)}  \,dt \Big)^{(2-\alpha)/2}\\
&\le& C \Big(\int_0^T |h_x|^{(2p-2)\alpha / (2-\alpha)}  \,dt \Big)^{(2-\alpha)/2}
\end{eqnarray*}
Fixing $\alpha=r/5$ yields the claim.
\end{proof}
%%
%%%%%%%%%%%%%%%%%%%%%%%%%%%%%%%%%%%%%%%%%%%%%%%%%%
\subsection{Blow up below criticality}

In this section we will study the blow up in a space below criticality, 
i.e.\ in some $H^s$ with $s<\frac12$.
This is a slight generalisation of Theorem \ref{t:regul}
and prepares the results of Leray-type shown later.

For $ \frac14\le \delta \le 1$ we obtain:
\begin{equation} \label{e:leray}
\begin{split}
 \tfrac12\partial_t |h|^2_{\delta}
&\leq -c  |h|^2_{{2+\delta}} + 2 \int_0^L (-\partial_x^2)^{\delta} h_x   \cdot  h_x h_{xx}\,dx\\
&\leq -c  |h|^2_{{2+\delta}} + C |h|_{{1+2\delta}}|h|_{{\frac94}}|h|_{{\frac54}}\\
&\leq -c  |h|^2_{{2+\delta}} + C |h|_{{2+\delta}}^{(9-2\delta)/4}  | h|_{\delta}^{(3+2\delta)/4}
\end{split}
\end{equation}
where we have used the Sobolev embedding $H^{\frac14}\subset L^4$.
\begin{remark}\label{r:interpol}
As it is used several times in the proofs,
we state the following elementary interpolation inequality. 
For $\gamma>\alpha$ and $\beta\in[\alpha,\gamma]$,
$$
|h|_{\beta} 
\le C |h|_{\alpha}^{\frac{\gamma-\beta}{\gamma-\alpha}}|h|_{\gamma}^{\frac{\beta-\alpha}{\gamma-\alpha}}.
$$
\end{remark}
Using interpolation between $H^\gamma$, $\gamma\le\tfrac54$ and $H^{2+\delta}$ implies
$$
 \tfrac12\partial_t |h|^2_{\delta}
 \le  -c  |h|^2_{{2+\delta}} + 
C | h|_{\gamma}^{(3+2\delta)/(4+2\delta-2\gamma)}  |h|_{{2+\delta}}^{(9+4\delta-6\gamma)/(4+2\delta-2\gamma)} 
$$
If we suppose $\gamma>\tfrac12$, then using Young inequality with $p=(8+4\delta-4\gamma)/(9+4\delta-6\gamma)$ and $q= (8+4\delta-4\gamma)/(2\gamma-1)$ 
we derive
$$
 \tfrac12\partial_t |h|^2_{\delta}
 \le 
C | h|_{\gamma}^{2(3+2\delta)/(2\gamma-1)} 
$$
We proved the following Theorem:
\begin{theorem}
Let $h\in C^\infty([0,t_0)\times[0,L])$ be 
a solution and fix $\gamma\in(\tfrac12,\tfrac54]$ and $\delta\in[\tfrac14,1]$.
Then  
$$ |h(t)|_{\delta} \to \infty\ \text{for}\  t\nearrow t_0
 \quad \Rightarrow \quad 
\int_0^{t_0} | h(t)|_{\gamma}^{2(3+2\delta)/(2\gamma-1)}\,dt =\infty\;.
$$ 
\end{theorem}
Note that for a blow up below criticality with $\delta<\frac12$ the 
$L^p([0,T],H^\gamma)$-norm in this theorem 
 has a smaller $p$ than assured by Theorem \ref{t:regul}.
The spaces in the above theorem should always have the 
same level of criticality.
%
%%
%%
%%%%%%%%%%%%%%%%%%%%%%%%%%%%%%%%%%%%%%%%%%%%%%%%%%%%%%%%%%%%%%%%%%%%%%%%%%%%%%%%%%%%%%%%%%%%%%%%%%%%%%%%%%%%%%%%%%%%
\section{Blow-up}
\label{sec:blowup}

In this section we discuss some properties of the blow up. 
First, at a possible  blow up time,
one expects that all norms with higher regularity 
than the critical norms will blow up,
in particular all $H^s$-norm with $s>1/2$ should blow up.
In Subsection \ref{ss:leray}, we give a lower bound on the blow-up 
in $H^s$-spaces, while in Subsection \ref{ss:singu}
we show a bound on the size of the set of singular times.
We illustrate that a blow up to $-\infty$ is 
more likely, but first we give some remarks on
possible shapes of a blow-up.
%%
%%%%%%%%%%%%%%%%%%%%%%%%%%%%%%%%%%%%%%%%%%%%%%%%%%%%%%%%%%%
\subsection{Some remarks}

Let us first give examples on which blow up profiles $v=h(\tau)$ are possible
at the blow up time $\tau$.
\begin{itemize}
\item If $v$ exhibits a {\em jump} like $\text{sign}(x)$,
then the Fourier-coefficients decay like $1/k$,
and thus $v$ is in $H^s$ if and only if $s<\frac12$.
\item If $v$ exhibits a {\em logarithmic pole} like $\log(|x|)$,
then the Fourier-coefficients decay like $1/k$,
and thus $v$ is in $H^s$ if and only if $s<\frac12$.
\item If $v$ exhibits a {\em cusp} like $|x|^\alpha$ for $\alpha\in(0,1)$,
then the Fourier-coefficients decay like $|k|^{-(1+\alpha)}$,
and thus $v$ is in $H^{1/2}$, and not a possible blow up.
\end{itemize}
%%
%%%%%%%%%%%%%%%%%%%%%%%%%%%%%%%%%%%%
\subsubsection{Stationary solutions}\label{sss:stationary}
The $L^2$ estimates~\eqref{e:energy} show that the only stationary solution
is $h \equiv 0$, as $|h(t)|_{L^2}\to0$ for $t\uparrow\infty$. On the other
hand the problem is one-dimensional, so it is worth trying to look for
solutions directly. The equation for stationary solutions is
$$
h_{xxxx}+(h_x^2)_{xx} = 0,
$$
so there are constants $A$, $B$ such that $h_{xx}+h_x^2=A x+B$. By the
periodic boundary conditions, $A=0$.

\emph{Case 1: $B=0$}. By direct computations, we get
$$
h(x) = c_1 + \log|1+c_2 x|,
$$
and the only periodic solution corresponds to $c_2=0$, a constant function.
Notice that, anyway, the solutions are singular with a log-like profile.

\emph{Case 2: $B=b^2$}. Again by direct computations,
$$
h(x) = c_1 + \log|\cosh b x + \frac{c_2}{b}\sinh b x|,
$$
and there are no periodic solutions. We remark that again the singularity
has a log-like profile.

\emph{Case 3: $B=-b^2$}. By elementary computations,
$$
h(x) = c_1 + \log|b\cos bx + c_2\sin bx|,
$$
all solutions are periodic on $[0,L]$ as long as $b=\frac{2\pi}{L}k$,
for some $k\in\mathbf{N}$. If $x_0$ is any zero of $b\cos bx+c_2\sin bx$,
we can write the solution as $h(x)=c_1 + \log|\sin b(x-x_0)|)$ (with a
different value of $c_1$). Again, the stationary profile is log-like.
%%
%%%%%%%%%%%%%%%%%%%%%%%%%%%%%%%%%%%%%%
\subsubsection{Self-similar solutions}
By exploiting the scaling \eqref{e:scaling}, we may look for solutions
of the following kind,
$$
h(t,x)=\varphi(\frac{x}{\root 4\of {T-t}}),
$$
where $\varphi$ is a suitable function. The equation for $h$ reads in terms
of $\varphi$ as
\begin{equation}\label{e:selfsimilar}
\varphi_{yyyy} + (\varphi_y^2)_{yy} + y\varphi_y = 0,
\qquad y\in\R,
\end{equation}
and, by the regularity of weak solutions 
one shows easily $\varphi$, $\varphi_{xx}\in L^2$ 
and hence $\varphi\in H^2(\R)$.
Here for simplicity we have neglected boundary conditions and formulated the problem
on the whole real line. The problem above can be recast in weak form as
$$
\int\varphi\eta_{yyyy}\,dy + \int\varphi_y^2\eta_{yy}\,dy - \int\varphi\eta\,dy - \int y\varphi\eta_y\,dy = 0,
\qquad\eta\in C_c^\infty,
$$
where the solution $\varphi\in H^1_\text{loc}(\R)$.

There is quite a strong numerical evidence that there are no solutions
to \eqref{e:selfsimilar} defined on the whole $\R$. This fact would rule
out self-similar solutions%
\footnote{Existence of self-similar solutions has been a long standing problem
for the Navier-Stokes equations. The problem was firstly posed by J. Leray
\cite{Ler34} in 1934 and finally solved by Ne\v{c}as,~R{\accent"17 u}\v{z}i\v{c}ka
\& \v{S}ver\'ak \cite{NecRuzSve96} in $1996$. Lately, Cannone \& Planchon~\cite{CanPla96}
proved existence of self-similar solution in Besov spaces.

Ne\v{c}as et al.~exploited a non-trivial maximum principle for $|u|^2+p$
(where $u$ is the velocity field and $p$ is the pressure). We remark
that no such fact seems to be true in this case.}.

%%
%%%%%%%%%%%%%%%%%%%%%%%%%%%%%%%%%%%%%%%%%%%%%%%%%%%%
\subsection{Leray-type results}
\label{ss:leray}
We will prove the following theorem, which is based on one of the several 
celebrated results of Leray~\cite{Ler34} on the Navier-Stokes equations. This relies mainly 
on a  comparison result for ODEs (see Lemma~\ref{l:techLR})
and energy estimates. It improves the results of Theorem \ref{t:regul},
which states that at blow-up for $s>\tfrac12$ the function
$t\to|h(t)|_s^{8/(2s-1)}$ is not integrable. The result now says
that it behaves like $\tfrac1t$. 
\begin{theorem}\label{t:leray}
Let $h\in C^\infty([0,t_0)\times [0,L])$ be a smooth local solution.
Then for $s>\frac12$
there is a universal constant $C>0$ such that
$|h(t)|_s \to \infty$ for $t\nearrow t_0$ (or for any subsequence)
implies  
$$|h(t)|_s \ge C(t_0-t)^{-(2s-1)/8}
\quad\text{for all } t\in[0,t_0)\;.
$$
\end{theorem}
\begin{proof}
We proceed by using energy estimates. 
Again use the notation $D=A^{1/4}=|\partial_x|$ and $B(u,v)=(u_x v_x)_{xx}$.

From (\ref{e:SG}) we obtain for $s=1+\delta$ with $\delta\in(-\tfrac12,\tfrac32)$
\begin{eqnarray*}
\partial_t|h|^2_{1+\delta}+ 2 |h|^2_{3+\delta} 
&=& -2 \int D^{2\delta}h_{xx} B(h,h)\,dx\\
&=& 4 \int D^{2\delta}h_{x} B(h,h_x)\,dx\\
&\le & C |h|_{1+\delta+\epsilon}|h|_{\tfrac52-\epsilon}|h|_{3+\delta}\;,
\end{eqnarray*}
where we used Proposition \ref{p:inequality} with
$\alpha=2+\delta$, $\beta=\tfrac12-\epsilon$, and $\gamma=-\alpha+\epsilon$
for some small $\epsilon\in(0,\tfrac12)$
such that $\epsilon+\delta \in (-\tfrac12,\tfrac32)$.
Now using interpolation (cf. Remark \ref{r:interpol}) yields
$$
\partial_t|h|^2_{1+\delta}+ 2 |h|^2_{3+\delta}
 \le C |h|_{1+\delta}^{\tfrac14(7-2\delta)} |h|_{3+\delta}^{\tfrac14(5+2\delta)}\;.
$$
As $(5+2\delta)<8$, we can apply Young's inequality 
with $p=8/(7-2\delta)$ and $q=8/(1+2\delta)$ 
to derive
$$ 
\partial_t |h|^2_{1+\delta}
\le C | h|_{1+\delta}^{2(5+2\delta)/(1+2\delta)}=C | h|_s^{2(3+2s)/(2s-1)}\;.
$$
Thus Lemma \ref{l:techLR} implies the theorem for $s\in(\tfrac12,\tfrac52)$.

Consider now $s=2+\delta$ with $\delta\in(-\tfrac12,\tfrac32)$.
\begin{eqnarray*}
\partial_t|h|^2_{2+\delta}+ 2 |h|^2_{4+\delta} 
&=& 2 \int D^{2\delta}h_{xxxx} B(h,h)\,dx\\
&=& -4 \int D^{2\delta}h_{xxx} B(h,h_x)\,dx\\
&=& -4 \int D^{2\delta}h_{xx} \left[B(h_x,h_x)+B(h,h_{xx})\right]\,dx\\
&\le & C |h|_{2+\delta+\epsilon}|h|_{\tfrac72-\epsilon}|h|_{3+\delta}+
       C |h|_{3+\delta+\epsilon}|h|_{\tfrac72-\epsilon}|h|_{2+\delta}\;,
\end{eqnarray*}
where we again used Proposition \ref{p:inequality}
with the same choice of $\alpha$, $\beta$, $\gamma$ and $\epsilon$.
Now using interpolation
$$
\partial_t|h|^2_{2+\delta}+ 2 |h|^2_{4+\delta}
 \le C |h|_{2+\delta}^{\tfrac14(7+2\delta)} |h|_{4+\delta}^{\tfrac14(5-2\delta)}\;
$$
and Young with $p=8/(5-2\delta)$ and $q=8/(3+2\delta)$ 
\begin{equation}\label{e:forsingular}
\partial_t |h|^2_{2+\delta}
\le C | h|_{2+\delta}^{2(7+2\delta)/(3+2\delta)}=C | h|_s^{2(3+2s)/(2s-1)}\;.
\end{equation}
Now Lemma \ref{l:techLR} finishes the proof for $s\in(\tfrac32,\tfrac72)$.

The general case is proven similarly, by distributing the derivatives
as evenly as possible on the trilinear terms, as in the proof of
Proposition~\ref{p:critimorereg}, and then applying Proposition~\ref{p:inequality},
possibly with different $\alpha$'s for different terms.
\end{proof}
\begin{remark}
\label{rem:leray}
We can also give a lower bound on the blow-up time $t_0$
depending on   $| h(0)|_{\delta}$ for $\delta>\frac12$.
To be more precise, using the upper bound in Lemma \ref{l:techLR}
the following is straightforward to verify. 
For all $s>\tfrac12$ there is a constant $c_s>0$ 
such that the solution is regular and smooth on $(t,t_\ast)$
if 
$ c_s |h(t)|_s^{8/(2s-1)}(t_\ast-t)<1.$

On the other hand, Theorem \ref{t:leray} immediately 
implies that near a blow up at $t_\ast$ we obtain 
for all $r\in(t,t_\ast)$, that 
$ c_s |h(r)|_s^{8/(2s-1)}(t_\ast-r) \ge 1.$
\end{remark}
%%
%%

%%
%%%%%%%%%%%%%%%%%%%%%%%%%%%%%%%%%%%%%%%%%%%%%%%%%%
\subsection{Criterion for point-wise blow up to \texorpdfstring{$-\infty$}{-infinity}}
\label{ss:negblow}

We show that for a blow up in $L^\infty$ the blow up to $-\infty$ 
is much more likely than the blow up to $\infty$.
This is mainly based on the a-priori estimate from Section \ref{sec:LF},
but first we use the following estimate:
$$
\begin{aligned}
\tfrac13 \partial_t \int_0^Lh^3\,dx
& = - \int_0^L h^2 h_{xxxx}\,dx - \int_0^L h^2 ((h_{x})^2)_{xx}\,dx\\
& =  2\int_0^L h h_x h_{xxx}\,dx + 4\int_0^L h (h_{x})^2 h_{xx}\,dx\\
& = - \int_0^L h (h_{xx})^2\,dx - \tfrac43 \int_0^L (h_{x})^4 \,dx,
\end{aligned}
$$
where we used the cancellation property~\eqref{e:key}.
Thus 
$$ \int_0^T \int_0^L (h_{x})^4\,dx\,dt 
\leq \int_0^L h^3(0)\,dx +  \int_0^T \int_0^L h^- (h_{xx})^2\,dx\,dt +  \int_0^L  h^- h^2\, dx.
$$
This implies:
\begin{theorem}
Let $h\in C^\infty([0,\tau)\times[0,L])$ be a smooth local solution.
If $\int_0^L h^3(0)\,dx$ is finite and $\| h\|_{L^4(0,\tau,W^{1,4})}=\infty$
then the negative part $h^-$ has to blow up. In other terms,
there are  $t_n\nearrow\tau$ and  $x_n\in[0,L]$ such that $  h(t_n,x_n) \to -\infty$.
\end{theorem}
\begin{corollary}
If $\int_0^L h^3(0)\,dx<\infty$ and $h^-$ uniformly bounded, then
$\| h\|_{L^4(0,T,W^{1,4})}<\infty$ and   $\int_0^T \int_0^L h^+ (h_{xx})^2\,dx\,dt<\infty$.
\end{corollary}
Let us now show that not only we have a point-wise
blow up, but also a blow up for some $\int_0^L e^{-\gamma h(t)}\,dx$,
while we know already by Section \ref{sec:LF}
that $\int_0^L e^{-\gamma h(t)}\,dx$ stays finite for $\gamma\in(0,2)$.
\begin{lemma}
Let $h\in C^\infty([0,\tau)\times[0,L])$ be a  smooth local solution.
If 
\begin{equation}
\label{e:exp_bl_ass}
 \int_0^T\int|h_x|^\alpha|h|^k\,dx\,dt \to\infty
\quad \text{for} \quad T\nearrow\tau
\end{equation}
 for some $\alpha\in(0,4)$ and $k\ge0$,
then 
$$ \int e^{-\gamma h(t)}\,dx  \to\infty 
\quad \text{for} \quad t\nearrow \tau
$$ 
for all $\gamma\in(0,2\alpha/(4-\alpha))$.
\end{lemma}

Note that the corresponding metric is always not critical.
It has less regularity.
Furthermore, note that for $\alpha\le 2$, by H\"older and interpolation,
the quantity in (\ref{e:exp_bl_ass}) will never blow up.

\begin{proof}
Using H\"older and results of Section \ref{sec:LF}
yields for any $\epsilon\in(0,\alpha/2)$ 
(i.e.\ $4\epsilon/\alpha\in(0,2)$),
$$
\begin{aligned}
\lefteqn{\Big(\int_0^T\int|h_x|^\alpha|h|^{k}\,dx\,dt \Big)^{\alpha/4}\leq}\qquad\\
&\leq  C \int_0^T \Big(\int|h_x|^\alpha|h|^{k}\,dx \Big)^{\alpha/4}\,dt\\
&\leq  C \int_0^T \int|h_x|^4 e^{4\epsilon h /\alpha}\,dx \cdot 
	\Big(\int e^{-4\epsilon h /(4-\alpha)} |h|^{4k/(4-\alpha)}\,dx \Big)^{(4-\alpha)/\alpha}\,dt\\
&\leq C\sup_{[0,T]} \Big( \int e^{-4\epsilon h /(4-\alpha)} |h|^{4k/(4-\alpha)}\,dx \Big)^{(4-\alpha)/\alpha}\\
&\leq C\sup_{[0,T]} \Big(\int e^{-\gamma h}\,dx \Big)^{(4-\alpha)/\alpha}
\end{aligned}
$$
for $\gamma\in (0, 4\epsilon  /(4-\alpha))$.
\end{proof}
%%
%%
%%
%%%%%%%%%%%%%%%%%%%%%%%%%%%%%%%%%%%%%%
\subsection{The set of singular times}
\label{ss:singu}

Let $h$ be a weak solution to~\eqref{e:SG} and consider the set of regular times of $h$,
$$
\regular = \{t\in(0,\infty): u\text{ is continuous with values in $H^1$ in a neighbourhood of }t\}.
$$
By Proposition~\ref{p:critimorereg}, $\regular$ is equal to the set of all times
$t$ such that $h$ is $C^\infty$ in space and time in a neighbourhood of $t$.
Define the set of singular times $\singular = [0,\infty)\setminus\regular$.

The next theorem proves (in the spirit of results of Leray~\cite{Ler34},
Scheffer~\cite{Sch76} for Navier-Stokes), that the set of singular times
is ``small''.
\begin{theorem}\label{t:singularset}
Given a weak solution $h$ to~\eqref{e:SG}, the set $\singular$ of singular times
of $h$ is a compact subset of $[0,\infty)$ and
$$
\mathcal{H}^\frac14(\singular) = 0,
$$
where $\mathcal{H}^\frac14$ is the $\tfrac14$-dimensional Hausdorff measure.
\end{theorem}
\begin{proof}
Fix a weak solution $h$ and define $\regular$ and $\singular$ as above.
The proof is divided in four steps.

\emph{1. $\singular$ is compact}.
The set $\regular$ is clearly open, hence $\singular$ is closed.
We prove that $\singular$ is bounded. Let $a_\bullet$ be the constant
given in part~\textbf{3} of Theorem~\ref{t:criticalexuniq}. Assume by
contradiction that $a_\bullet<|h(t)|_{\frac12}$ for all $t\geq0$.
By interpolation and using the energy inequality~\eqref{e:energy},
$$
a_\bullet^{\frac83} t
< \int_0^t |h(s)|_{\frac12}^{\frac83}\,ds
\leq \int_0^t |h(s)|_{L^2}^{\frac23} |h(s)|_{2}^2\,ds
\leq |h(0)|_{L^2}^{\frac23}\int_0^t |h(s)|_2^2\,ds
\leq 2 |h(0)|_{L^2}^{\frac83}.
$$
Hence for some $t_0>0$, $|h(t_0)|_{\frac12}\leq a_\bullet$ and Theorems
\ref{t:criticalexuniq} and \ref{t:weakstrong} imply that the solution
$h$ is regular in $[t_0,\infty)$.

\emph{2. $\singular$ has Lebesgue measure $0$}.
As any open set of $\R$ is the countable union of
disjoint open intervals we have $\regular=\bigcup_j I_j$, 
where the open intervals $I_j$ are the connected
components of $\regular$.

Define $\singular_2=\{t:u(t)\not\in\dot H^2\}$.
Trivially, $\regular\subset\singular_2^c$, hence $\singular_2\subset\singular$.
If $t_0\in\singular\setminus\singular_2$, by Proposition~\ref{p:critimorereg}
$t_0$ is the endpoint of some $I_j$, hence $\singular\setminus\singular_2$ is
at most countable. Finally, the energy estimate \eqref{e:energy} 
implies that $\singular_2$ has measure $0$.

\emph{3. Estimate on the length of bounded $I_j$}.
Indeed, let $I_j$ be a bounded component of $\regular$ and let $t_1,t_2\in I_j$.
From Remark \ref{rem:leray} we know  $c(t_2-s)|h(s)|_2^{8/3}\geq1$,
and hence $c(t_2-s)^{-3/4}\leq |h(s)|_2^2$, for $s\in(t_1,t_2)$.
Integrating for $s\in(t_1,t_2)$ and using the energy inequality~\eqref{e:energy},
yields
\begin{equation}\label{e:intervalbound}
c(t_2-t_1)^{\frac14}
\leq \int_{t_1}^{t_2} |h(s)|_2^2\,ds
\leq \frac12|h(0)|_{L^2}^2.
\end{equation}

\emph{4. $\mathcal{H}^{\frac14}(\singular)=0$}.
Write $I_j=(a_j,b_j)$ for bounded intervals. From~\eqref{e:intervalbound}
it follows that
$$
\sum_j (b_j - a_j)^{\frac14}<\infty,
$$
while $\sum_j (b_j - a_j)<\infty$, by the first step of the proof. Now we
can proceed as in the proof of Theorem 2 of~\cite{Sch76} to get the conclusion.
\end{proof}
%%
%%
%%
%%%%%%%%%%%%%%%%%%%%%%%%%%%%%%%%%%%%%%%%%%%%%%%%%%%%%%%%%%%%%%%%%%%%%%%%%%%%%%%%%%%%%%%%%%%%%%%%%%%%%%%%%%%%%%%%%%%%
\appendix
\section{An inequality for the non-linearity}
Given three real numbers $\alpha$, $\beta$, $\gamma$, consider the following condition.
\begin{condition}\label{e:condinequality}
The real numbers $\alpha$, $\beta$, $\gamma$ satisfy
\begin{itemize}
\item $\alpha$, $\beta\geq0$,
\item $\alpha+\beta+\gamma\geq\frac12$ with strict inequality if at least one is equal to $\frac12$,
\item if $\gamma<0$, then either at least one of $\alpha$ and $\beta\leq\frac12$, or at least one $\geq-\gamma$.
\end{itemize}
\end{condition}
\begin{lemma}\label{l:sommafacile}
For every $\gamma\in\R$ there is $c>0$ such that for every $a\geq1$,
$$
\sum_{|k|\leq a}|k|^{-2\gamma}
\leq\begin{cases}
ca^{1-2\gamma} &\qquad \gamma<\frac12,\\
c\log a &\qquad \gamma=\frac12,\\
c & \qquad\gamma>\frac12.
\end{cases}
$$
\end{lemma}
\begin{lemma}\label{l:sommadifficile}
Let $\alpha$, $\beta$ and $\gamma$ satisfy \eqref{e:condinequality} (with $\alpha\leq\frac12$ or $\beta\geq-\gamma$ when $\gamma<0$). Then there is $c>0$ such that for each $m\in\Z$, with $m\neq0$,
$$
\sum_{\substack{|k|<2|m|\\0<|k-m|<\frac12|k|}}\frac1{|k-m|^{2\alpha}|k|^{2\gamma}}\leq c|m|^{2\beta}.
$$
\end{lemma}
\begin{proof}
Notice that, if $|k|<2|m|$ and $|k-m|<\frac12|k|$, then $\frac23|m|\leq|k|<2|m|$, since
$$
\frac23|m|\leq\frac23|k-m|+\frac23|k|<\frac13|k|+\frac23|k|=|k|.
$$
Then apply Lemma \ref{l:sommafacile}.
\end{proof}
Consider
$$
B(u,v)=(u_xv_x)_{xx}.
$$
\begin{proposition}\label{p:inequality}
If $\alpha$, $\beta$ and $\gamma$ satisfy \eqref{e:condinequality},
there exists $c>0$ such that for all $u\in \dot H^{1+\alpha}$,
$v\in \dot H^{1+\beta}$ and $w\in\dot H^{2+\gamma}$,
$$
\langle B(u,v), w\rangle
\leq c|u|_{1+\alpha}|v|_{1+\beta}|w|_{2+\gamma}.
$$ 
\end{proposition}
\begin{proof}
\emph{Step 1.}
Write the functions $u$, $v$, $w$ in the Fourier expansion,
$$
u=\sum_{k\neq0} u_k \e^{ikx}
$$
(and similarly for $v$ and $w$), so that
$$
(u_xv_x)_{xx}=\sum_{k\neq0} k^2\Bigl(\sum_{l+m=k}lmu_lv_m\Bigr)\e^{ikx}
$$
and by Cauchy-Schwartz,
$$
\begin{aligned}
\langle B(u,v), w\rangle
&=    \sum_{k\neq0} k^2\overline{w_k}\Bigl(\sum_{l+m=k}lmu_lv_m\Bigr)\\
&\leq |w|_{2+\gamma}\Bigl[\sum_{k\neq0} |k|^{-2\gamma}\Bigl(\sum_{l+m=k} |l m u_l v_m|\Bigr)^2\Bigr]^{\frac12}.
\end{aligned}
$$
Hence, it is sufficient to analyse only the second term in the above product. Set for every $k\neq0$,
\begin{gather*}
A_k=\{(l,m):l+m=k,\ |l|\geq\frac12|k|,\ |m|\geq\frac12|k|\},\\
B_k=\{(l,m):l+m=k,\ |l|<\frac12|k|\},\\
C_k=\{(l,m):l+m=k,\ |m|<\frac12|k|\},
\end{gather*}
and, for simplicity, $U_l=|l|^{1+\alpha}|u_l|$ and $V_m=|m|^{1+\beta}|v_m|$.

\emph{Step 2.} We start by analysing the sum on $A_k$.
\begin{align*}
\sum_{k\neq0} |k|^{-2\gamma}\Bigl(\sum_{A_k}|l m u_l v_m|\Bigr)^2
&\leq\sum_{k\neq0} |k|^{-2\gamma}\Bigl(\sum_{A_k}|l|^{-\alpha}|m|^{-\beta}U_lV_m\Bigr)^2\\
	\textit{\footnotesize (using Young's inequality)\quad}
&\leq c\sum_{k\neq0} |k|^{-2\gamma}\Bigl(\sum_{A_k}|l|^{-\alpha-\beta}U_lV_m\Bigr)^2\\
&\quad +c\sum_{k\neq0} |k|^{-2\gamma}\Bigl(\sum_{A_k}|m|^{-\alpha-\beta}U_lV_m\Bigr)^2,\\
\intertext{the two terms are similar, we bound only the first one,}
	\textit{\footnotesize(using Cauchy inequality)\quad}
&\leq c\sum_{k\neq0} |k|^{-2\gamma}\Bigl(\sum_{A_k}|l|^{-2(\alpha-\beta)}U_l^2\Bigr)\Bigl(\sum_{A_k}V_m^2\Bigr)\\
	\textit{\footnotesize(switching the sums)\quad}
&\leq c|v|_{1+\beta}^2\sum_{l\neq0}|l|^{-2(\alpha+\beta)}U_l^2\Bigl(\sum_{|k|\leq 2|l|}|k|^{-2\gamma}\Bigr)\\
	\textit{\footnotesize(using Lemma \ref{l:sommafacile})\quad}
&\leq c|u|_{1+\alpha}^2|v|_{1+\beta}^2.
\end{align*}
\emph{Step 3.} Next, we analyse the sum on $B_k$ (the sum on $C_k$ being entirely similar).
Notice that, when using Cauchy inequality below, we are free to weigh either the terms in $u$
or in $v$ with derivatives. We shall choose one of the two depending on the values of $\gamma$
(wherever we need an exponent to be $\leq\frac12$ or $\geq-\gamma$, according to condition
\eqref{e:condinequality}).
\begin{align*}
\sum_{k\neq0} |k|^{-2\gamma}\Bigl(\sum_{B_k}|l m u_l v_m|\Bigr)^2
&\leq\sum_{k\neq0} |k|^{-2\gamma}\Bigl(\sum_{B_k}|l|^{-\alpha}|m|^{-\beta}U_lV_m\Bigr)^2\\
	\textit{\footnotesize(using Cauchy inequality)\quad}
&\leq\sum_{k\neq0} |k|^{-2\gamma}\Bigl(\sum_{B_k}|l|^{-2\alpha}|m|^{-2\beta}V_m^2\Bigr)\Bigl(\sum U_l^2\Bigr)\\
	\textit{\footnotesize(switching the sums)\quad}
&\leq|u|_{1+\alpha}^2\sum_{m\neq0}|m|^{-2\beta}V_m^2\Bigl(\sum_{\substack{|k|<2|m|\\ 0<|k-m|<\frac12|k|}}|k-m|^{-2\alpha}|k|^{-2\gamma}\Bigr)\\
	\textit{\footnotesize(using Lemma \ref{l:sommadifficile})\quad}
&\leq c|u|_{1+\alpha}^2|v|_{1+\beta}^2.
\end{align*}
The proof is complete.
\end{proof}

%%%%%%%%%%%%%%%%%%%%%%%%%%%%%%%%%%%%%%%%%%%%%%%%%%%%%%%%%%%%%%%%%%%%%%%%%%%%%%%%%%%%%%%%

\section{Blow up for ODEs}

%%%%%%%%%%%%%%%%%%%%%%%%%%%%%%%%%%%%%%%%%%%%%%%%%%%%%%%%%%%%%%%%%%%%%%%%%%%%%%%%%%%%%%%%

The following elementary lemma is crucial to prove Leray-type bounds.
We state and proof it for completeness.
\begin{lemma}\label{l:techLR}
Let $\varphi:(0,t_0)\to\mathbb{R}$ be a non-negative  function such that for $p>1$  we have
$ \partial_t \varphi \le C\varphi^p$, on $ (0, t_0)$.

Then,
$ \varphi(t_n) \uparrow \infty $ for a subsequence $t_n\uparrow t_0$,
implies 
$$\varphi(t) \ge [(p-1)C(t_0-t)]^{-1/(p-1)}
\quad\text{for all }t\in(0,t_0)\;.$$

Moreover,
$$
\varphi(t)
\le \Big[\varphi(s)^{-(p-1)}+C(p-1)s - C(p-1)t\Big]^{-1/(p-1)}
\quad \text{for all } 0<s<t<t_0\;.
$$
\end{lemma}
\begin{proof}
We have for $t>s$
$$
\tfrac1{p-1}(\varphi(s)^{-(p-1)}- \varphi(t)^{-(p-1)})
=\int_{\varphi(s)}^{\varphi(t)} \frac1{z^p} dz 
= \int_s^t \frac{\partial_t \varphi}{\varphi^p}d\tau 
\le C(t-s)
$$
Now for $t_n\uparrow t_0$ we obtain
$$\tfrac1{p-1}\varphi(t)^{-(p-1)} \le C(t_0-t)$$
and finally
$$\varphi(t) \ge [(p-1)C(t_0-t)]^{-1/(p-1)}$$
for all $t\in(0,t_0)$.

For the second result
$$
\varphi(s)^{-(p-1)}- C(p-1)(t-s)
\le \varphi(t)^{-(p-1)}
$$
and thus 
$$
\varphi(t)
\le \Big[\varphi(s)^{-(p-1)}+C(p-1)s - C(p-1)t\Big]^{-1/(p-1)}\;.
$$
\end{proof}
%%
%%
%%
%%%%%%%%%%%%%%%%%%%%%%%%%%%%%%%%%%%%%%%%%%%%%%%%%%%%%%%%%%%%%%%%%%%%%%%%%%%%%%%%%%%%%%%%
%%%%%%%%%%%%%%%%%%%%%%%%%%%%%%%%%%%%%%%%%%%%%%%%%%%%%%%%%%%%%%%%%%%%%%%%%%%%%%%%%%%%%%%%
\section{Analytic Semigroups}

The following properties of analytic semigroups are well known,
but we give short sketches of proofs for the sake of completeness.
\begin{lemma}\label{lem:SG}
Consider $A=\partial_x^4$ subject to periodic boundary conditions on $[0,L]$
and $T>0$. For all $u=\sum_k u_k e_k\in L^2$ and $\alpha>0$,
$$
|s^\alpha A^\alpha\e^{-sA} u|_{L^2}\to 0,
\qquad\text{for }s\to0,
$$
and for all $f\in L^\infty(0,T,L^2)$, with $0$ average on $(0,1)$,
and $1 + a - b=0$, with $a>-1$ and $b<1$,
$$
I_{a,b}(f)(t) = \int_0^ts^a A^b e^{-(t-s)A} f(s)\,ds
$$
converges to $0$ in $L^2$ as $t\to0$.
\end{lemma}
\begin{proof}
The first statement is obvious by Lebesgue theorem, since
$$
|s^\alpha A^\alpha \e^{-sA} u|_{L^2}^2
\leq C \sum_k (k^4s)^{2\alpha} e^{-csk^4} u_k^2.
$$
For the second statement note that (with a change of variables)
$|I_{a,b}(f)(t)|_{L^2}^2$ is equal to
$$
\begin{aligned}
\lefteqn{\sum_{k\neq0}k^{8b}\int_0^t\int_0^t s^a r^a k^{8b} \e^{-c(2t-s-r)k^4} f_k(s)f_k(r) d\,ds\,dr = }\\
& = t^{2a+2}\int_0^1\int_0^1 (1-s)^a (1-r)^a \sum_{k\not=0} k^{8b} \e^{-c(s + r)tk^4} f_k(t-ts )f_k(t-tr ) \,ds\,dr\\
& = t^{2b}\int_0^1\int_0^1 (1-s)^a(1-r)^a \Bigl[\Bigl(\sum_{(s+r)tk^4\leq\epsilon}+\sum_{(s+r)tk^4>\epsilon}\Bigr) k^{8b}\e^{-c(s+r)tk^4}f_k(t-ts)f_k(t-tr)\Bigr]\,ds \,dr \\
& \leq c\int_0^1\int_0^1\frac{(1-s)^a (1-r)^a}{(s+r)^{2b}}\Bigl(\epsilon^{2b}\|f\|_{L^\infty(L^2)} + \sum_{(s +r )tk^4>\epsilon}f_k(t-ts )f_k(t-tr )\Bigr)\,ds\,dr,
\end{aligned}
$$
which goes to zero by Lebesgue theorem if one first takes the limit as $t\to0$ and then as $\epsilon\downarrow0$,
since the function $(1-r)^a(1-s)^a(r+s)^{-2b}$ is integrable and the other term is bounded for $\epsilon\leq1$
and $t\leq T$.
\end{proof}
%%
%%
%%
%%%%%%%%%%%%%%%%%%%%%%%%%%%%%%%%%%%%%%%%%%%%%%%%%%%%%%%%%%%%%%%%%%%%%%%%%%%%%%%%%%%%%%%%%%%%%%%%%%%%%%%%%%%%%%%%%%%

%%
%%
%%
%%%%%%%%%%%%%%%%%%%%%%%%%%%%%%%%%%%%%%%%%%%%%%%%%%%%%%%%%%%%%%%%%%%%%%%%%%%%%%%%%%%%%%%%%%%%%%%%%%%%%%%%%%%%%%%%%%%
\end{document}